\documentclass[11pt, a4paper]{article}

\usepackage[utf8]{inputenc}
\usepackage[T1]{fontenc}
\usepackage{geometry}
\usepackage{amsmath, amssymb, amsthm, amsfonts}
\usepackage{mathtools}
\usepackage{graphicx}
\usepackage{float}
\usepackage{hyperref}
\hypersetup{hidelinks}
\usepackage{enumitem}
\usepackage{booktabs}
\usepackage{algorithm}
\usepackage{algpseudocode}
\usepackage{bm}
\usepackage{xcolor,array}
\usepackage{subcaption}
\graphicspath{{figures/}}
\title{\textbf{Shape-Derivative Moment Approximation for a Reduced Time-Harmonic Biot Poroelastic System with Random Interfaces}}
\author{Yongle Hao\thanks{Email: haoyl@zknu.edu.cn}\\ 
	\small School of Mathematics and Statistics, Zhoukou Normal University\\
	\small Zhoukou, HeNan, China\\\\
	Fuchang Huo\thanks{
		Email: fuchang.huo@polyu.edu.hk; li.jz@sustech.edu.cn. Corresponding authors.}\\ 
	\small Department of Applied Mathematics, The Hong Kong Polytechnic University\\
	\small Kowloon, Hong Kong, China\\\\
	Jingzhi Li\footnotemark[2]\\
	\small Department of Mathematics, Southern University of Science and Technology\\
	\small Shenzhen, Guangdong, China\\\\
	Kai Zhang\thanks{Email: zhangkaimath@jlu.edu.cn}\\
	\small Department of Mathematics, Jilin University\\
	\small Changchun, Jilin, China
}
\date{November 2026}

\numberwithin{equation}{section}

\newtheorem{theorem}{Theorem}[section]
\newtheorem{lemma}[theorem]{Lemma}

\newtheorem{assumption}[theorem]{Assumption}
\newtheorem{remark}[theorem]{Remark}

\newcommand{\bfu}{\mathbf{u}}

\newcommand{\bfn}{\mathbf{n}}
\newcommand{\bfx}{\mathbf{x}}
\newcommand{\bfy}{\mathbf{y}}
\newcommand{\dd}{\mathrm{d}}

\newcommand{\Cor}{\operatorname{Cor}}

\newcommand{\E}{\mathbb{E}}
\newcommand{\Prob}{\mathbb{P}}
\newcommand{\PT}{\mathbf{P}_\Gamma}

\newcommand{\jump}[1]{[#1]}

\begin{document}
	\maketitle
	
	\begin{abstract}
	We develop first- and second-moment approximations for a reduced time-harmonic Biot poroelastic interface system with random material interfaces. The model retains mixture inertia, quasi-static Darcy discharge, and explicit storage and coupling coefficients. The interface separates two media, each with constant Lam\'{e} parameters and hydraulic mobility. We use shape calculus to analyze random perturbations normal to the interface. We derive the first shape derivative in a total-traction and normal-flux formulation, establish a Taylor expansion on a fixed function space under a non-resonance assumption, and prove error bounds on fixed observation regions. The mean approximation is second-order accurate, and the covariance and variance approximations are third-order accurate. We estimate the variance error directly. Applying pivoted Cholesky compression to the covariance matrix of the coefficients in a fixed basis of smooth interface modes reduces the computation of the leading covariance to deterministic sensitivity solves. Numerical experiments in two and three dimensions examine perturbation rates, low-rank truncation, coupled statistics, and computational cost.

	\end{abstract}

	

	\section{Introduction}\label{sec:introduction}
	
	Wave propagation in fluid-saturated porous media is governed by coupled deformation--flow mechanisms.  In Biot's theory, the elastic displacement of the porous skeleton is coupled to a pore-pressure field through the pressure gradient in the momentum balance and the volumetric strain in the mass balance \cite{Biot1956}.  In applications such as seismic imaging, acoustic liners, biological tissues, and engineered porous composites, material interfaces are rarely known exactly.  Across such interfaces, the elastic moduli and permeability may change abruptly, while the interface location itself may be affected by measurement error, manufacturing tolerance, or unresolved microscopic geometry.  The resulting forward model is therefore a coupled interface problem with geometric uncertainty.
	
	The numerical treatment of random domains and random interfaces is more delicate than that of random coefficients on a fixed domain.  Mapping each realization to a nominal configuration generally produces transformed coefficients that depend nonlinearly on the random geometry and may have reduced spatial regularity \cite{Canuto2007}.  Classical Monte Carlo methods are robust but require many deterministic solves, whereas stochastic Galerkin and collocation methods can be highly efficient only when the parametric dependence is sufficiently regular and the effective stochastic dimension is moderate \cite{Babuska2004,Babuska2007,Ghanem1991,Nobile2008,Xiu2006}.  Recent work on random domains, moving domains, and domain uncertainty has clarified how regularity of the domain map controls the approximation of statistics for elliptic, parabolic, elastic, and wave models \cite{Djurdjevac2021,HarbrechtKarnaevSchmidlin2024,HenriquezLabarca2026}.  In parallel, modern stochastic Galerkin and low-rank strategies have continued to improve the practical solution of high-dimensional stochastic PDE discretizations \cite{BachmayrEigelEisenmannVoulis2025,BachmayrVoulis2022,DuswaldKeithLazarovPetridesWohlmuth2024,Harbrecht2012,Youngnoi2025}.
	
	For small geometric perturbations, shape-calculus perturbation methods offer a complementary route \cite{Delfour2011,Sokolowski1992}.  The random solution is expanded around the deterministic solution on a nominal interface, and leading-order statistics are computed from deterministic sensitivity equations. This idea is closely related to earlier shape-differentiability techniques for inverse obstacle scattering and discontinuity identification \cite{Hettlich1995,Hettlich1998}, and has been developed further for elliptic, electromagnetic, and scattering interface problems \cite{BaoCaoHaoZhang2018,BaoMaLaiLi2026,Dolz2020,Hao2018,HarbrechtLi2013,HiptmairLi2018,HiptmairLi2013}. Recent work on shape Taylor expansions for wave scattering provides a closely related higher-order expansion perspective for deterministic scattering problems \cite{BaoMaLaiLi2026}.
	The present work develops a corresponding moment analysis for a {reduced low-frequency} time-harmonic Biot-type interface system { retaining mixture inertia and quasi-static Darcy discharge \cite{Osinov2023,ZienkiewiczShiomi1984}}. This extension is not a direct scalar generalization. The state contains both a vector displacement and a scalar pressure, while the interface laws couple displacement, pressure, total traction, and Darcy flux. Consequently, second-moment statistics are also intrinsically coupled: the covariance structure contains displacement-displacement, displacement-pressure, and pressure-pressure correlations simultaneously.

	{This paper presents} a shape-calculus and covariance-based moment analysis for a time-harmonic Biot-type interface system with random normal perturbations. Across the material interface, the Lam\'{e} parameters and the permeability are assumed to be piecewise constant.
	{ Specifically, the main contributions are:}
	\begin{enumerate}[label=(\roman*),leftmargin=*,itemsep=0.4em,topsep=0.5em]
		\item {\textbf{Coupled interface sensitivities.}}
		We derive the first shape derivative in broken Sobolev spaces and formulate the resulting interface conditions in total-traction and flux form. This formulation keeps the sign convention, the variation of the normal vector, and the role of possible source jumps explicit, and {  is well suited to numerical implementation. The weak formulation uses a fixed test space and lifts only the displacement and pressure jumps.}
		
		\item { \textbf{Moment estimates for complex-valued fields.}}
		Under { the stated non-resonance and regularity assumptions, and for centered, uniformly bounded perturbations},  we obtain asymptotic approximations for the first two statistical moments: the nominal solution approximates the mean with an $O(\epsilon^2)$ error, while { $\epsilon^2$ times} the covariance of the first shape derivative approximates the solution covariance with an $O(\epsilon^3)$ error. { These are absolute error bounds on fixed physical observation regions separated from the interface.} A key feature of the formulation is the coupled covariance structure induced by the Biot system. {The Hermitian covariance kernel includes the mixed displacement--pressure blocks, and we estimate the variance error directly using $H^1\hookrightarrow L^4$.}
		
		\item { \textbf{Low-rank finite element realization.}
			We combine pivoted Cholesky compression of a fixed smooth-mode coefficient covariance with sensitivities obtained by differentiating the interface-fitted finite element equations. The sensitivity solves reuse the nominal matrix factorization. The error analysis separates state and sensitivity approximation, geometric consistency, and covariance truncation; experiments in two and three dimensions examine these effects and the role of mixed displacement--pressure statistics.}
	\end{enumerate}
	{ The first-order second-moment framework and pivoted Cholesky compression are established methods \cite{HarbrechtLi2013,Hao2018,Harbrecht2012}; the contribution of this work is a physically consistent realization of these methods for the coupled Biot interface problem, which constitutes a more complex system.}

	Recent numerical analysis of Biot and poromechanics systems provides the discretization context for this work, including mixed, discontinuous Galerkin, stabilized, and splitting schemes for poroelastic (and coupled) models \cite{AntoniettiBonettiBotti2023,AntoniettiBottiMazzieriNatiPoltri2022,BoonFumagalliScotti2023,GuoSunWangYueZheng2025,KreuzerZanotti2026,PeDeLaRivaGasparHuAdlerRodrigoZikatanov2025,GomezVargasMardalRuizBaierVinje2023}.  The present work takes a different approach: we examine how a random perturbation of the interface propagates into the first and second moments of the reduced Biot solution. The numerical section therefore relies on standard deterministic solvers solely as a tool for testing the perturbation expansion and the low-rank covariance approximation. 
	
	The paper is organized as follows.  Section~\ref{section:Model} introduces the deterministic and stochastic reduced Biot interface model.  Section~\ref{section:Shape_Calculus} derives the first shape derivative.  Section~\ref{section:Statistics} establishes {  the mean, covariance, and variance approximations.}  Section~\ref{section:Numerical_Implementation} describes the low-rank numerical realization.  Section~\ref{section:Error_Estimates} records the finite element and truncation error structure.  Section~\ref{section:Numerical_experiments} presents numerical experiments, and Section~\ref{section:Conclusion} offers concluding remarks.
	
	\section{Model problem}\label{section:Model}
	
	\subsection{Reduced deterministic Biot interface system}

\begin{figure}[htbp]
		\centering
		\includegraphics[width=0.48\textwidth]{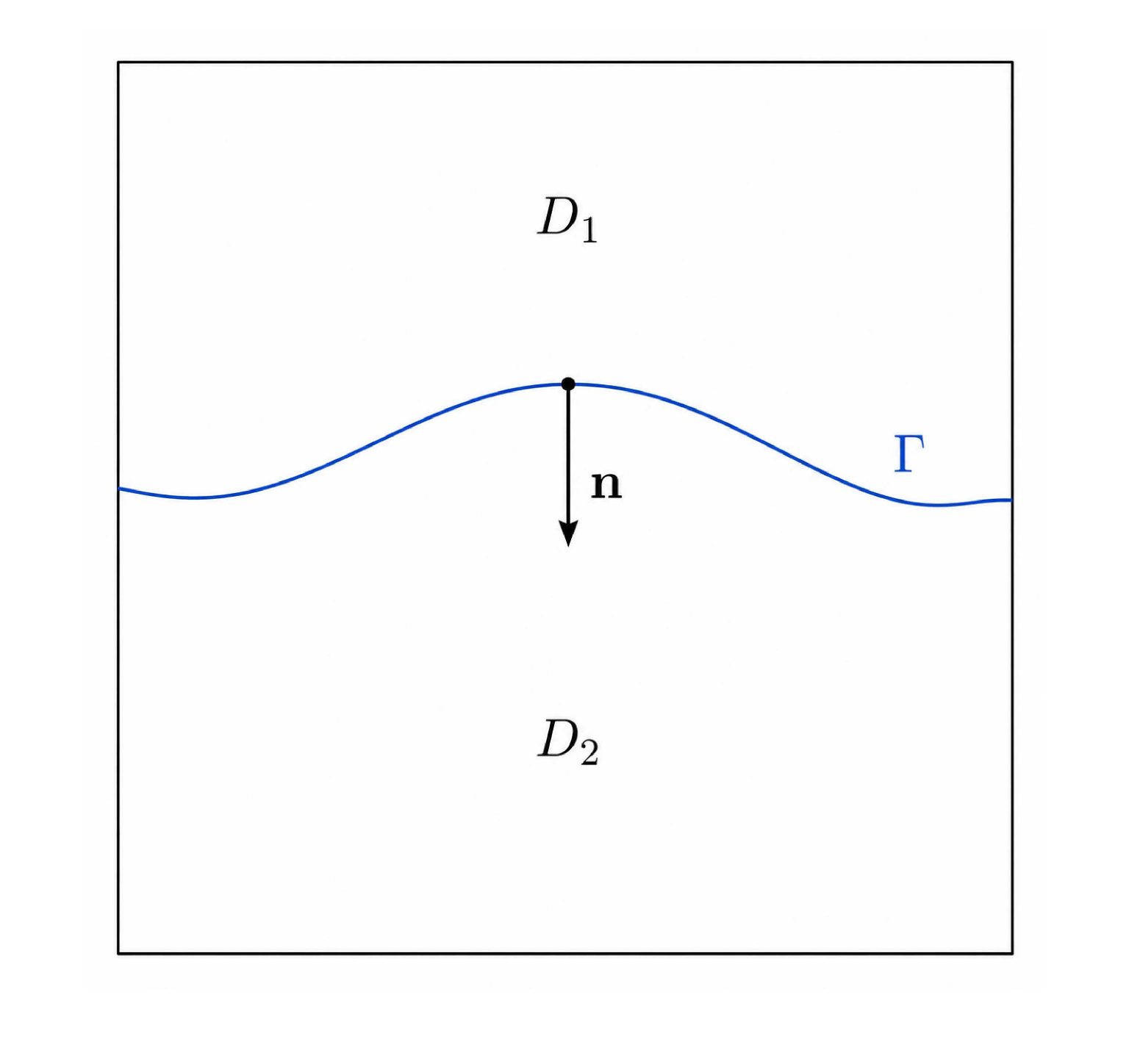}
		\caption{Geometry of the deterministic interface $\Gamma$.}
		\label{fig:reference_interface}
	\end{figure}

	Consider a bounded Lipschitz domain $D\subset\mathbb{R}^d$ ($d=2$ or $3$). Let $\Gamma$ denote the nominal interface, which divides $D$ into two open subdomains $D_1$ and $D_2$ with $\overline D = \overline D_1 \cup \overline D_2$ and $\Gamma = \partial D_1 \cap \partial D_2$. { Assume that $D_1$ and $D_2$ are Lipschitz subdomains. The moving part of their common interface satisfies the stronger regularity and support conditions stated in Section~\ref{sec:stochastic_geometry}.} The unit normal $\bfn$ is fixed to point from $D_1$ to $D_2$, and for a quantity $w$ with traces $w_i=w|_{D_i}$ we use
	\[
	\jump{w}=w_2-w_1 .
	\]
	When piecewise Sobolev regularity is required, we use the notation
	\[
	H^s(D_1\cup D_2):=\{v:\ v_i\in H^s(D_i),\ i=1,2\}.
	\]
	The reference geometry is illustrated in Figure~\ref{fig:reference_interface}.
	
	The model retains the inertia of the fluid--solid mixture but treats the Darcy discharge quasi-statically, neglecting both the relative fluid--solid acceleration term in the momentum balance and the inertial correction to the Darcy discharge. These assumptions give the approximation that Osinov calls UP2 \cite{Osinov2023}. For the underlying dynamic theory and the development of reduced formulations, see Biot \cite{Biot1956} and Zienkiewicz and Shiomi \cite{ZienkiewiczShiomi1984}, respectively. The interface analysis presented below does not quantify the error introduced by this low-frequency approximation.

	Before nondimensionalization, the governing equations are given by
	\begin{align}
		\rho_{\rm phys}\partial_{tt}\mathbf u_{\rm phys}
		-\nabla\cdot\bigl(\boldsymbol\sigma_{\rm phys}(\mathbf u_{\rm phys})
		-\alpha p_{\rm phys}\mathbf I\bigr)
		&=\mathbf f_{\rm phys},
		\label{eq:physical_time_momentum}\\
		c_{\rm phys}\partial_t p_{\rm phys}
		+\alpha\nabla\cdot\partial_t\mathbf u_{\rm phys}
		-\nabla\cdot(K_{\rm phys}\nabla p_{\rm phys})
		&=g_{\rm phys}.
		\label{eq:physical_time_mass}
	\end{align}
	Here, $\rho_{\rm phys}$ denotes the mixture density, $K_{\rm phys}$ is the hydraulic mobility (the intrinsic permeability divided by the dynamic viscosity of the fluid), and $-K_{\rm phys}\nabla p_{\rm phys}$ is the Darcy discharge. The Biot--Willis coefficient $\alpha$ and the storage coefficient $c_{\rm phys}$ are retained explicitly.

	Let $L_0$, $E_0$, and $T_0$ be positive reference scales for length, stress, and time, respectively. We introduce the dimensionless variables
	\[
	\mathbf x=\frac{\mathbf x_{\rm phys}}{L_0},
	\quad
	t=\frac{t_{\rm phys}}{T_0},
	\quad
	\mathbf u=\frac{\mathbf u_{\rm phys}}{L_0},
	\quad
	p=\frac{p_{\rm phys}}{E_0}.
	\]
	The corresponding dimensionless material parameters, angular frequency, and source terms are defined by
	\begin{align}\label{eq:nondimensionalization}
		\lambda_i&=\lambda_{{\rm phys},i}/E_0, \hspace{-20mm} &
		\mu_i&=\mu_{{\rm phys},i}/E_0, \hspace{-20mm}&
		c_i&=c_{{\rm phys},i}E_0,\notag\\ 
		\kappa_i&=\frac{K_{{\rm phys},i}E_0T_0}{L_0^2},\hspace{-20mm} &
		\rho_i&=\frac{\rho_{{\rm phys},i}L_0^2}{E_0T_0^2},\hspace{-20mm} &
		\theta&=\theta_{\rm phys}T_0,\\
		\mathbf f&=\mathbf f_{\rm phys}L_0/E_0,\hspace{-20mm} &
		g&=g_{\rm phys}T_0. \hspace{-20mm} & \notag
	\end{align}
	The Biot--Willis coefficient $\alpha_i$ is dimensionless and remains unchanged. The same reference scales are used in all material subdomains, ensuring that the physical transmission conditions are preserved under nondimensionalization.

	Adopting the time-harmonic convention $e^{i\theta t}$, we obtain the following system for the complex amplitudes $U=(\mathbf u,p)$:
	\begin{equation}\label{eq:biot_pde}
		\begin{cases}
			-\nabla\cdot\bigl(\boldsymbol\sigma_i(\mathbf u)
			-\alpha_i p\mathbf I\bigr)
			-\theta^2\rho_i\mathbf u
			=\mathbf f
			&\text{in }D_i,\\[2pt]
			-\nabla\cdot(\kappa_i\nabla p)
			+i\theta\bigl(c_i p+\alpha_i\nabla\cdot\mathbf u\bigr)
			=g
			&\text{in }D_i,
		\end{cases}
	\end{equation}
	where the elastic stress tensor and linear strain tensor are given by
	\begin{align*}
		\boldsymbol\sigma_i(\mathbf u) =2\mu_i\boldsymbol\varepsilon(\mathbf u) +\lambda_i(\nabla\cdot\mathbf u)\mathbf I, \quad \boldsymbol\varepsilon(\mathbf u) =\frac12\bigl(\nabla\mathbf u+\nabla\mathbf u^T\bigr).
	\end{align*}
	The material coefficients are real and piecewise constant, satisfying
	\begin{align*}
		\mu_i>0,\quad d\lambda_i+2\mu_i>0,\quad \kappa_i>0,\quad \rho_i>0,\quad c_i>0,\quad 0<\alpha_i\le1.
	\end{align*}
	Throughout the analysis, the material parameters and the angular frequency $\theta>0$ are fixed. No parameter-uniform estimates are claimed in the nearly incompressible, impermeable, or resonant limits.

	We assume perfect mechanical bonding and hydraulic contact across the material interface $\Gamma$, together with homogeneous clamped and drained conditions on the exterior boundary $\partial D$. The corresponding transmission and boundary conditions are
	\begin{equation}\label{eq:interface_cond}
		\begin{gathered}
			\jump{\mathbf u}=0,\quad
			\jump{p}=0,\quad
			\jump{\mathcal A(U)\mathbf n}=0,\quad
			\jump{\mathcal Q(U)\cdot\mathbf n}=0
			\quad\text{on }\Gamma,\\
			\mathbf u=0,\quad p=0
			\quad\text{on }\partial D,
		\end{gathered}
	\end{equation}
	where
	\begin{equation}\label{eq:total_stress_flux}
		\mathcal A_i(U)
		=\boldsymbol\sigma_i(\mathbf u)-\alpha_i p\mathbf I,
		\quad
		\mathcal Q_i(U)=\kappa_i\nabla p.
	\end{equation}
	Here, $\mathcal A_i(U)$ denotes the total stress, while $-\mathcal Q_i(U)$ represents the Darcy discharge. The transmission conditions enforce continuity of displacement, pressure, total traction, and normal Darcy discharge. In particular, continuity of pressure and total traction yields
	\[
	\jump{\boldsymbol\sigma(\mathbf u)\mathbf n}
	=\jump{\alpha}p\mathbf n.
	\]
	Thus, a discontinuity in the Biot--Willis coefficient generally leads to a jump in the elastic traction, even though the total traction remains continuous. This distinction plays an essential role in the subsequent shape differentiation of the transmission conditions.

	\subsection{Interface conditions and weak formulation}

	We formulate the weak problem in the complex Hilbert space
	\[
	X=H^1_0(D;\mathbb C^d)\times H^1_0(D;\mathbb C),
	\quad
	\|U\|_X^2
	=\|\mathbf u\|_{H^1(D)}^2
	+\|p\|_{H^1(D)}^2.
	\]
	Let $X^*$ denote the continuous antilinear dual of $X$. For $U=(\mathbf u,p)\in X$ and $W=(\mathbf v,q)\in X$, we define the sesquilinear form
	\begin{align}
		b(U,W)
		=\sum_{i=1}^2\int_{D_i}\bigl\{&
		2\mu_i\boldsymbol\varepsilon(\mathbf u):
		\overline{\boldsymbol\varepsilon(\mathbf v)}
		+\lambda_i(\nabla\cdot\mathbf u)
		\overline{\nabla\cdot\mathbf v}
		-\theta^2\rho_i\mathbf u\cdot\overline{\mathbf v}
		-\alpha_i p\,\overline{\nabla\cdot\mathbf v}
		\notag\\
		&+\kappa_i\nabla p\cdot\overline{\nabla q}
		+i\theta c_i p\overline q
		+i\theta\alpha_i(\nabla\cdot\mathbf u)\overline q
		\bigr\}\,d\mathbf x,
		\label{eq:bilinear_form}
	\end{align}
	and the antilinear functional
	\[
	\ell(W)
	=\int_D\bigl(
	\mathbf f\cdot\overline{\mathbf v}
	+g\overline q
	\bigr)\,d\mathbf x.
	\]
	The nominal weak problem is to find $U_0\in X$ such that
	\begin{equation}\label{eq:weak_state}
		b(U_0,W)=\ell(W),
		\quad\forall\,W\in X.
	\end{equation}
	The continuity of displacement and pressure across $\Gamma$ is incorporated into the global space $X$, while the continuity of total traction and normal Darcy flux is imposed weakly through the variational formulation.

	\begin{assumption}[Fixed data and non-resonance]
		\label{ass:data}
		The source terms $\mathbf f$ and $g$ are deterministic, fixed in physical coordinates, and admit $C^2$ extensions to a neighborhood of $\overline D$. Moreover, the nominal weak operator is injective in the sense that
		\[
		b(U,W)=0\quad\forall\,W\in X
		\quad\Longrightarrow\quad U=0.
		\]
	\end{assumption}
	
	The non-resonance assumption excludes frequencies at which the nominal coupled operator has a nontrivial kernel. Such an assumption cannot be justified solely by Darcy dissipation, since divergence-free displacement modes with vanishing pressure are not directly damped by the Darcy term. Lemma~\ref{lem:fredholm} establishes the Fredholm property of the coupled operator, and Theorem~\ref{thm:pullback} shows that invertibility persists under sufficiently small interface perturbations.

	\subsection{Stochastic interface perturbation}\label{sec:stochastic_geometry}
	
\begin{figure}[htbp]
		\centering
		\includegraphics[width=0.48\textwidth]{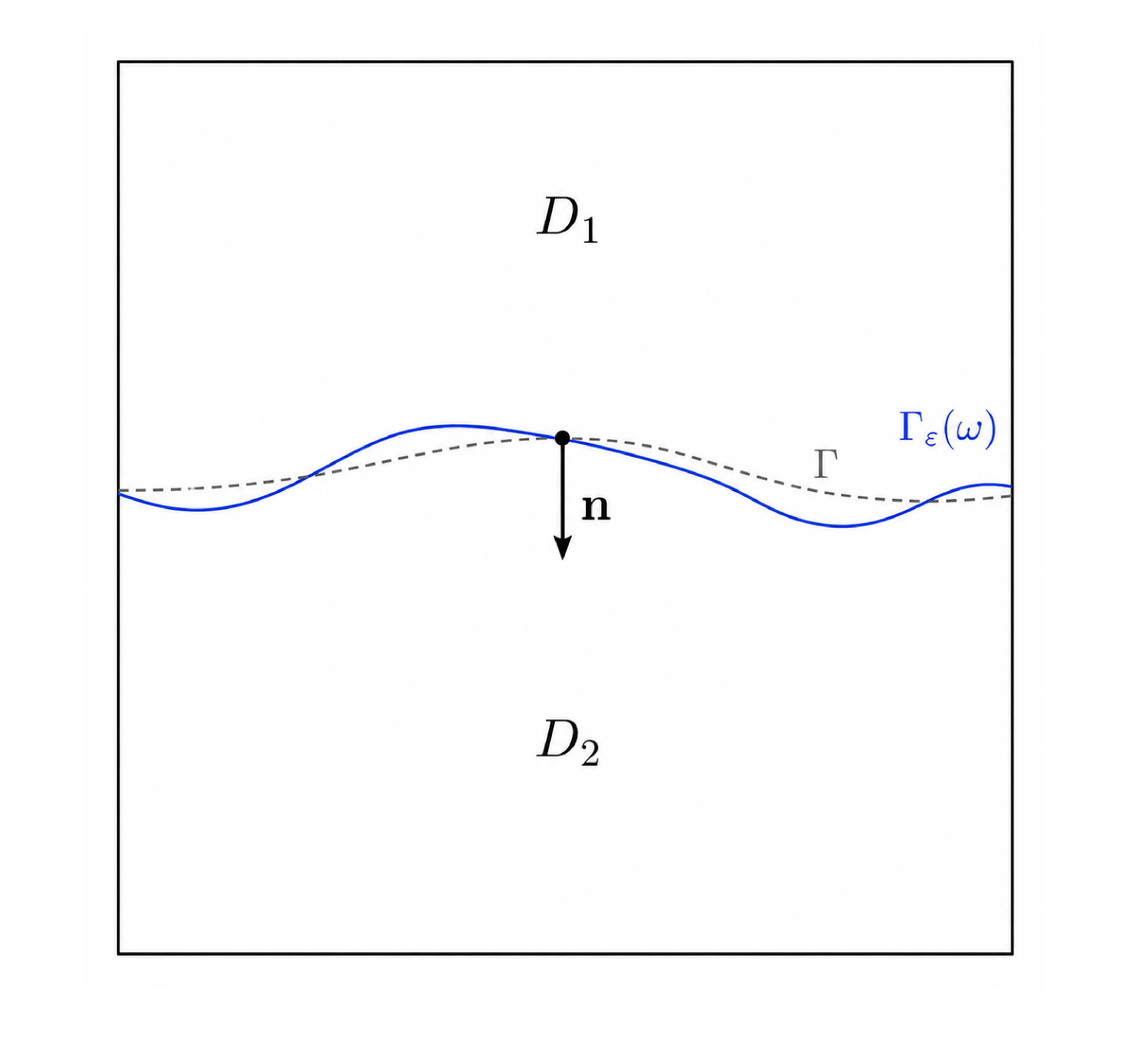}
		\caption{Geometry of the random interface perturbation $\Gamma_\epsilon(\omega)$.}
		\label{fig:random_interface}
	\end{figure}

	Let $(\Xi,\mathcal F,\Prob)$ be a probability space, and let $\omega\in\Xi$ denote a random outcome.  The random interface is a normal graph over $\Gamma$,
	\begin{equation*}
		\Gamma_\epsilon(\omega)=\{\bfx+\epsilon\beta(\bfx,\omega)\bfn(\bfx):\ \bfx\in\Gamma\},
	\end{equation*}
	as illustrated in Figure~\ref{fig:random_interface}. Here $0<\epsilon\ll1$ is the perturbation amplitude.  We assume
	\begin{equation*}
		\E[\beta](\bfx)=0,
		\quad
		\Cor_\beta(\bfx,\bfy)=\E[\beta(\bfx,\omega)\beta(\bfy,\omega)].
	\end{equation*}

	The active portion of $\Gamma$ is assumed to lie in an embedded $C^{3,1}$ surface patch (or a finite union of such patches). All admissible perturbation amplitudes are supported in a fixed compact subset of this smooth portion contained in $D$. For a smooth closed interface contained in $D$, this subset may be all of $\Gamma$. In particular, if $\Gamma$ intersects $\partial D$, the amplitudes vanish in a common neighborhood of $\Gamma\cap\partial D$. Let $B$ be the closed real subspace of $C^{2,1}$ amplitudes supported in this fixed compact set, equipped with the $C^{2,1}$ norm defined in local surface charts. These amplitudes are extended by zero to the inactive portion of $\Gamma$. We assume that $\beta$ is strongly measurable as a $B$-valued random variable, with its expectation understood in the Bochner sense in $B$. Expectations of the response fields are taken in the respective spaces specified below.

	For the pathwise shape calculus, we use the stronger regularity assumption
	\begin{equation} \label{eq:beta_assumption}
		\beta\in L^\infty(\Xi;B),
		\quad
		\|\beta\|_{L^\infty(\Xi;B)}\le M_\beta.
	\end{equation}
	This ensures a realization-independent tubular neighborhood and a uniform diffeomorphism radius for small $\epsilon$.  The weaker assumption $\beta\in L^2(\Xi;B)$ is sufficient for some moment identities but not for uniform $O(\epsilon^2)$ remainders. When $\Gamma$ has boundary on $\partial D$, all admissible perturbations are understood to satisfy the fixed-end compatibility condition stated above; in numerical examples this is enforced either by a smooth endpoint cut-off or by using a closed interface.
	
	Fix a Lipschitz observation set $K$, possibly a finite disjoint union of Lipschitz sets, such that $\overline K\Subset D\setminus\Gamma$. In particular, $K$ is separated from the nominal interface by a non-zero distance.
	Choose a fixed bounded real-linear extension operator
	\[
	E:B\longrightarrow W^{2,\infty}(D;\mathbb R^d),\quad
	V_\beta=E\beta,\quad V_\beta|_\Gamma=\beta\mathbf n,
	\quad \|V_\beta\|_{W^{2,\infty}(D)}\le C_E\|\beta\|_B.
	\]
	The extension is supported in a fixed tubular neighborhood of the moving interface and vanishes in fixed neighborhoods of $\partial D$ and $\overline K$. Such an extension is obtained in tubular coordinates by extending the amplitude and normal along normal rays and applying a fixed smooth cutoff. For an interface meeting the exterior boundary, the fixed support condition on $B$ permits the same construction away from the intersection. The extension and its norm bound are independent of the realization.
	
	Set
	\begin{equation}\label{eq:geometry_map}
		F_{\epsilon,\beta}(\mathbf x)=\mathbf x+\epsilon V_\beta(\mathbf x),
		\quad D_{i,\epsilon}=F_{\epsilon,\beta}(D_i),\quad
		\Gamma_\epsilon=F_{\epsilon,\beta}(\Gamma).
	\end{equation}
	There exists a geometric radius $\epsilon_{\rm geo}>0$ such that $F_{\epsilon,\beta}$ is an orientation-preserving bi-Lipschitz diffeomorphism for every $\beta\in B$ with $\|\beta\|_B\le M_\beta$ and every $|\epsilon|<\epsilon_{\rm geo}$. The material coefficients are transported with their respective subdomains, whereas the physical source terms remain fixed in the spatial coordinates. For $|\epsilon|\le\epsilon_0$, with $0<\epsilon_0<\epsilon_{\rm geo}$ as in Theorem~\ref{thm:pullback}, let $U_\epsilon$ denote the perturbed physical solution, and define its pullback by
	\begin{equation}\label{eq:pullback_state}
		\widetilde U_\epsilon=U_\epsilon\circ F_{\epsilon,\beta}\in X.
	\end{equation}

	The principal statistical errors are evaluated on the fixed observation set $K$ defined above. By the choice of extension, $V_\beta$ vanishes in a fixed neighborhood of $\overline K$, uniformly in $\beta$. Consequently,
	\begin{equation}\label{eq:fixed_observation_identity}
		F_{\epsilon,\beta}|_K=\mathrm{id},\quad
		\widetilde U_\epsilon|_K=U_\epsilon|_K.
	\end{equation}
	Thus, solutions are compared at identical physical observation coordinates for every admissible perturbation and hence for almost every $\omega\in\Xi$ by~\eqref{eq:beta_assumption}. In contrast, whole-domain pullback statistics involve the material derivative, as distinguished below.

	\section{Shape calculus for the reduced Biot system}\label{section:Shape_Calculus}

	Let $V_\beta$ and $F_{\epsilon,\beta}$ be as defined in Section~\ref{sec:stochastic_geometry}. The material and Eulerian shape derivatives are defined piecewise in the two materials by
	\[
	\dot U[\beta]=\left.\frac{d}{d\epsilon}(U_\epsilon\circ F_{\epsilon,\beta})\right|_{\epsilon=0},
	\quad U'[\beta]=\dot U[\beta]-\nabla U_0 V_\beta.
	\]
	We write $U':=(\dd\mathbf u,\dd p)=(\mathbf u',p')$. Appendix~\ref{app:fixed_domain}proves that the material derivative $\dot U[\beta]$ exists in the fixed global space $X$, together with a uniform Taylor remainder. For the strong interface identities, we further assume that $U_0$ is piecewise $H^3$ near the velocity support. Let $\PT=\mathbf I-\mathbf n\otimes\mathbf n$, extend $\mathbf n$ constantly along normal rays, and interpret surface divergence row-wise. The normal variation is
	\begin{equation}\label{eq:normal_variation}
		\mathrm d\mathbf n=-\nabla_\Gamma\beta.
	\end{equation}

	The following two local identities are the basic Hadamard formulas used below. They isolate the geometric part of the shape derivative from the particular structure of the Biot operator.
	
	\begin{lemma}\label{lem:hadamard_traction}

		Let $A_{\epsilon,i}$ be tensor fields defined on $D_{i,\epsilon}$. Suppose that their one-sided pullbacks are differentiable at $\epsilon=0$ in the local $C^0$ norm up to the interface, and that the nominal fields $A_i=A_{0,i}$ are $C^1$ up to the interface on each side. Let $A'$ denote the phasewise Eulerian derivative. If $[A_\epsilon\mathbf n_\epsilon]=0$ on $\Gamma_\epsilon$ for all sufficiently small $\epsilon$ and $V|_\Gamma=\beta\mathbf n$, then

		\begin{equation}\label{eq:hadamard_traction}
			\jump{A'\bfn}
			=
			\nabla_\Gamma\cdot\bigl(\beta\jump{A\PT}\bigr)
			-
			\beta\jump{\nabla\cdot A},
		\end{equation}
		where $\nabla_\Gamma\cdot$ is applied row-wise.
	\end{lemma}
	
	\begin{proof}
		The transported conormal condition reads $\jump{A_\epsilon\bfn_\epsilon}(x+\epsilon\beta\bfn)=0$. Differentiating at $\epsilon=0$ and writing the Eulerian derivative rather than the material derivative gives
		\begin{equation}\label{eq:conormal_derivative_step}
			\jump{A'\bfn}=-\beta\jump{\partial_n(A\bfn)}-\jump{A\dd\bfn}.
		\end{equation}
		Since $\dd\bfn=-\nabla_\Gamma\beta$ is tangential, $A\dd\bfn=A\PT\dd\bfn$, and hence
		\[
		-\jump{A\dd\bfn}=\jump{A\PT}\nabla_\Gamma\beta .
		\]
		For a tensor field $A$, the row-wise surface decomposition is
		\begin{equation}\label{eq:tensor_surface_divergence}
			\nabla\cdot A=\partial_n(A\bfn)+H A\bfn+\nabla_\Gamma\cdot(A\PT),
		\end{equation}
		where $H=\nabla_\Gamma\cdot\bfn$. Taking the jump in \eqref{eq:tensor_surface_divergence} and using $\jump{A\bfn}=0$ gives
		\[
		\jump{\partial_n(A\bfn)}=\jump{\nabla\cdot A}-\nabla_\Gamma\cdot\jump{A\PT}.
		\]
		Substitution into \eqref{eq:conormal_derivative_step} yields
		\[
		\jump{A'\bfn}
		= -\beta\jump{\nabla\cdot A}+\beta\nabla_\Gamma\cdot\jump{A\PT}+\jump{A\PT}\nabla_\Gamma\beta
		=\nabla_\Gamma\cdot(\beta\jump{A\PT})-\beta\jump{\nabla\cdot A},
		\]
		which proves \eqref{eq:hadamard_traction}.
	\end{proof}
	
	\begin{lemma}\label{lem:hadamard_flux}

		Suppose that the one-sided pullbacks of $p_{\epsilon,i}$ are differentiable at $\epsilon=0$ in the local $C^1$ norm up to the interface, and that $p=p_0$ is $C^2$ up to $\Gamma$ on each side, with $[p]=0$ on $\Gamma$. Assume that $\kappa_i$ is constant in each moving subdomain $D_{i,\epsilon}$. If $[\kappa\nabla p_\epsilon\cdot\mathbf n_\epsilon]=0$ on $\Gamma_\epsilon$ for all sufficiently small $\epsilon$ and $V|_\Gamma=\beta\mathbf n$, then the Eulerian derivative $\dd p$ satisfies

		\begin{equation}\label{eq:hadamard_flux}
			\jump{\kappa\nabla \dd p\cdot\bfn}
			=
			\nabla_\Gamma\cdot\bigl(\beta\jump{\kappa}\nabla_\Gamma p\bigr)
			-
			\beta\jump{\nabla\cdot(\kappa\nabla p)} .
		\end{equation}
	\end{lemma}

	\begin{proof}
		{  Let $r_\epsilon=\kappa_i\nabla p_{\epsilon,i}$ in $D_{i,\epsilon}$.} Differentiating the transported flux condition $\jump{r_\epsilon\cdot\bfn_\epsilon}(x+\epsilon\beta\bfn)=0$ gives
		\[
		\jump{\kappa\nabla\dd p\cdot\bfn}
		=-\beta\jump{\partial_n(\kappa\partial_n p)}-\jump{\kappa\nabla p\cdot\dd\bfn}.
		\]
		Using \eqref{eq:normal_variation} and the continuity of $p$ along $\Gamma$, the tangential gradient $\nabla_\Gamma p$ is single-valued and
		\[
		-\jump{\kappa\nabla p\cdot\dd\bfn}=\jump{\kappa}\nabla_\Gamma p\cdot\nabla_\Gamma\beta .
		\]
		The scalar surface decomposition
		\[
		\nabla\cdot(\kappa\nabla p)
		=\partial_n(\kappa\partial_n p)+H\kappa\partial_n p+\nabla_\Gamma\cdot(\kappa\nabla_\Gamma p)
		\]
		and the flux continuity $\jump{\kappa\partial_n p}=0$ imply
		\[
		\jump{\partial_n(\kappa\partial_n p)}
		=
		\jump{\nabla\cdot(\kappa\nabla p)}-
		\nabla_\Gamma\cdot(\jump{\kappa}\nabla_\Gamma p).
		\]

		Substituting the preceding identity into the differentiated flux condition gives
		\begin{align*}
			[\kappa\nabla \dd p\cdot \bfn]
			={}&-\beta[\partial_n(\kappa\partial_n p)]
			+[\kappa]\nabla_\Gamma p\cdot\nabla_\Gamma\beta\\
			={}&-\beta[\nabla\cdot(\kappa\nabla p)]
			+\beta\nabla_\Gamma\cdot
			([\kappa]\nabla_\Gamma p)
			+[\kappa]\nabla_\Gamma p\cdot\nabla_\Gamma\beta\\
			={}&-\beta[\nabla\cdot(\kappa\nabla p)]
			+\nabla_\Gamma\cdot
			(\beta[\kappa]\nabla_\Gamma p),
		\end{align*}
		where the last equality follows from the surface product rule. This proves \eqref{eq:hadamard_flux}.

	\end{proof}

	The Eulerian shape derivative and the material derivative are related by
	\begin{equation}\label{eq:shape_material_relation}
		U'[\beta]=\dot U[\beta]-\nabla U_0\,V_\beta
		\quad\text{piecewise in }D_1\cup D_2.
	\end{equation}
	The local piecewise $H^3$ assumption stated above does not imply global $H^3$ regularity at exterior corners.

	\begin{lemma}\label{lemma:shape_deriv}
		Fix $\beta\in B$ and suppose that Assumption~\ref{ass:data} and the geometric assumptions of Section~\ref{sec:stochastic_geometry} hold. Assume further that $U_0$ is piecewise $H^3$ in a neighborhood of $\operatorname{supp}V_\beta$. Then the Eulerian shape derivative $U'=(\mathrm d\mathbf u,\mathrm dp) =\dot U[\beta]-(\nabla U_0)V_\beta$ belongs to the broken $H^1$ space and satisfies the homogeneous bulk equations associated with \eqref{eq:biot_pde} in each $D_i$, with $\mathrm d\mathbf u=\mathbf{0}$ and $\mathrm dp=0$ on $\partial D$. The following transmission conditions hold in the trace sense, with the traction and flux terms on the left-hand sides understood as weak conormal traces:
		\begin{align}
			[\mathrm d\mathbf u]
			&= -\beta[\partial_n\mathbf u_0],
			\label{eq:sd_u}\\
			[\mathrm dp]
			&= -\beta[\partial_n p_0],
			\label{eq:sd_p}\\
			[\boldsymbol\sigma(\mathrm d\mathbf u)\mathbf n
			-\alpha\,\mathrm dp\,\mathbf n]
			&= \nabla_\Gamma\cdot
			\bigl(\beta[\mathcal A(U_0)\mathbf P_\Gamma]\bigr)
			+\beta[\mathbf f+\theta^2\rho\mathbf u_0],
			\label{eq:sd_traction_total}\\
			[\kappa\nabla\mathrm dp\cdot\mathbf n]
			&= \nabla_\Gamma\cdot
			\bigl(\beta[\kappa]\nabla_\Gamma p_0\bigr)
			+\beta[g]
			-i\theta\beta[c p_0+\alpha\nabla\cdot\mathbf u_0].
			\label{eq:sd_flux_general}
		\end{align}
		Here $\mathcal A_i(U_0)=\boldsymbol\sigma_i(\mathbf u_0)
		-\alpha_i p_0\mathbf I$,
		$\mathbf P_\Gamma=\mathbf I-\mathbf n\otimes\mathbf n$,
		and tensor surface divergence is understood row-wise.
		Assumption~\ref{ass:data} implies $[\mathbf f]=[g]=0$.
	\end{lemma}
	
	\begin{proof}
		Appendix~\ref{app:fixed_domain} establishes the differentiability of the pulled-back solution in $X$. Together with the assumed local regularity, this yields $U'=\dot U[\beta]-(\nabla U_0)V_\beta$ in the broken $H^1$ space. The convergence of the Eulerian difference quotient in $L^2_{\mathrm{loc}}(D_i)$ allows passage to the limit in the physical bulk equations, yielding the homogeneous equations in each $D_i$, since the material coefficients and physical sources are fixed. The homogeneous exterior conditions follow from the fixed boundary data and the vanishing of $V_\beta$ near $\partial D$.
		
		Differentiating the transported continuity conditions and using $V_\beta=\beta\mathbf n$ on $\Gamma$ gives
		\[
		[\mathrm d\mathbf u+\beta\partial_n\mathbf u_0]=0, \quad [\mathrm dp+\beta\partial_n p_0]=0,
		\]
		which proves \eqref{eq:sd_u}--\eqref{eq:sd_p}.
		
		To differentiate the conormal transmission conditions, we use the weak normal-trace formulation. Write $F_\epsilon=F_{\epsilon,\beta}$, $G_\epsilon=(DF_\epsilon)^{-1}$, and $J_\epsilon=\det DF_\epsilon$. Applying the Piola transform to each row of the stress tensor and to the flux vector gives
		\[
		\widehat{\mathcal A}_\epsilon
		=J_\epsilon(\mathcal A(U_\epsilon)\circ F_\epsilon)G_\epsilon^T,
		\quad
		\widehat{\mathcal Q}_\epsilon
		=J_\epsilon G_\epsilon(\mathcal Q(U_\epsilon)\circ F_\epsilon).
		\]
		Their $L^2(D_i)$ differentiability follows from the fixed-space $H^1$ differentiability. By the Piola identities and the assumed regularity of the source terms, their divergences are also differentiable in $L^2(D_i)$. Hence the transforms are differentiable in $H(\operatorname{div};D_i)$, and the continuity of the normal-trace operator permits differentiation of the transmission conditions in the weak trace sense. The assumed local regularity and Lemmas~\ref{lem:hadamard_traction} and~\ref{lem:hadamard_flux} yield the stated interface expressions.
		
		By Lemma~\ref{lem:hadamard_traction},
		\[
		[\boldsymbol\sigma(\mathrm d\mathbf u)\mathbf n
		-\alpha\,\mathrm dp\,\mathbf n]
		=
		\nabla_\Gamma\cdot
		\bigl(\beta[\mathcal A(U_0)\mathbf P_\Gamma]\bigr)
		-\beta[\nabla\cdot\mathcal A(U_0)].
		\]
		Using the bulk equation
		\[
		\nabla\cdot\mathcal A_i(U_0)
		=-\mathbf f-\theta^2\rho_i\mathbf u_{0,i}
		\]
		yields \eqref{eq:sd_traction_total}.
		
		Similarly, Lemma~\ref{lem:hadamard_flux} gives
		\[
		[\kappa\nabla\mathrm dp\cdot\mathbf n]
		=
		\nabla_\Gamma\cdot
		\bigl(\beta[\kappa]\nabla_\Gamma p_0\bigr)
		-\beta[\nabla\cdot(\kappa\nabla p_0)].
		\]
		Substituting
		\[
		\nabla\cdot(\kappa_i\nabla p_{0,i})
		=i\theta(c_i p_{0,i}
		+\alpha_i\nabla\cdot\mathbf u_{0,i})-g
		\]
		completes the proof of \eqref{eq:sd_flux_general}.
	\end{proof}

	Let $X_{\rm br}$ denote the broken $H^1$ product space with zero exterior trace, equipped with the square root of the sum of the squared $H^1$ norms of its components over the two subdomains, and define $b_{\rm br}:X_{\rm br}\times X\to\mathbb C$ by the same volume integrals as $b$. Let $Y_\Gamma$ be the trace of $H_0^1(D)$ on $\Gamma$ with its quotient norm, and write $H^{-1/2}(\Gamma):=Y_\Gamma^*$ for its continuous antilinear dual. The duality pairing is antilinear in the test trace. For a smooth closed interior interface, $Y_\Gamma=H^{1/2}(\Gamma)$ with equivalent norms. This definition also covers interfaces meeting $\partial D$. Vector traces and their duals are understood componentwise.

	\begin{theorem}	\label{thm:sd_wellposed}
		Suppose Assumption~\ref{ass:data} holds. Let $g_u\in H^{1/2}(\Gamma;\mathbb C^d)$, $g_p\in H^{1/2}(\Gamma;\mathbb C)$, $g_T\in H^{-1/2}(\Gamma;\mathbb C^d)$, and $g_J\in H^{-1/2}(\Gamma;\mathbb C)$. If $\Gamma$ meets $\partial D$, assume that all four jump data are supported in a fixed compact set $\Gamma_*\Subset\Gamma\setminus\partial D$. The constants below may depend on $\Gamma_*$ and its distance from $\partial D$.
		Then there exists a unique $U'\in X_{\rm br}$ satisfying the homogeneous bulk equations and the prescribed interface jumps in the weak sense, with
		\begin{equation}\label{eq:transmission_stability}
			\|U'\|_{X_{\rm br}}\leq C\bigl(
			\|g_u\|_{H^{1/2}}+\|g_p\|_{H^{1/2}}
			+\|g_T\|_{H^{-1/2}}+\|g_J\|_{H^{-1/2}}\bigr).
		\end{equation}
		In particular, the jump data in Lemma~\ref{lemma:shape_deriv} uniquely determine $U'[\beta]$ for each admissible direction $\beta$. Moreover, if the map from $\beta\in B$ to the jump data is bounded and real-linear, then so is $\beta\mapsto U'[\beta]$.
	\end{theorem}
	
	\begin{proof}
		The trace theorem on the Lipschitz subdomains provides bounded right inverses. Hence, there exists a lifting $L=(L_u,L_p)\in X_{\rm br}$ such that
		$[L_u]=g_u$, $[L_p]=g_p$, and
		\[
		\|L\|_{X_{\rm br}}\leq
		C\bigl(\|g_u\|_{H^{1/2}}+\|g_p\|_{H^{1/2}}\bigr).
		\]
		Such a lifting can be constructed on one side of the interface, with zero extension on the other. If $\Gamma$ meets $\partial D$, a cut-off supported near $\Gamma_*$ ensures zero exterior trace. Only the primal jumps are lifted.
		
		Define $Z\in X$ by
		\begin{equation}\label{eq:weak_jump_lifting}
			b(Z,W)=-\langle g_T,\mathbf v|_\Gamma\rangle
			-\langle g_J,q|_\Gamma\rangle-b_{\rm br}(L,W),
			\quad W=(\mathbf v,q)\in X.
		\end{equation}
		The right-hand side is bounded, and Lemma~\ref{lem:fredholm} yields a unique $Z\in X$. Setting $U'=Z+L$ gives the prescribed primal jumps.
		
		Testing with functions supported in each material subdomain yields the homogeneous bulk equations:
		\[
		\nabla\cdot\mathcal A(U')=-\theta^2\rho\mathbf u'\in L^2,
		\quad
		\nabla\cdot\mathcal Q(U')
		=i\theta(c p'+\alpha\nabla\cdot\mathbf u')\in L^2.
		\]
		Hence the stress rows and flux belong to $H(\operatorname{div};D_i)$, with well-defined normal traces on $\partial D_i$. Pairing these traces with restrictions of global test functions defines interface functionals on $Y_\Gamma$. Green's identity and \eqref{eq:weak_jump_lifting} then yield the prescribed conormal jumps $g_T$ and $g_J$. The negative signs follow from $\mathbf n_1=\mathbf n$, $\mathbf n_2=-\mathbf n$, and $[w]=w_2-w_1$.
		
		The bounded inverse and trace estimates yield \eqref{eq:transmission_stability}. Uniqueness follows since the difference of two solutions belongs to $X$ and solves the nominal homogeneous problem.
		
		For the final assertion, the local piecewise $H^3$ hypothesis in Lemma~\ref{lemma:shape_deriv} is understood for every admissible direction $\beta\in B$; no common $H^3$ neighborhood or uniform local $H^3$ norm is needed. Define the real-linear map
		\[
		M:B\longrightarrow X_{\rm br},\quad
		M\beta=(\nabla U_0)E\beta.
		\]
		This map is well defined: the stated local regularity controls the product near $\operatorname{supp}E\beta$, while $E\beta$ vanishes away from its support and in a fixed neighborhood of the exterior boundary. To show that $M$ has a closed graph, suppose $\beta_n\to\beta$ in $B$ and $M\beta_n\to Z$ in $X_{\rm br}$. Since $E:B\to L^\infty(D;\mathbb R^d)$ is bounded and $\nabla U_0\in L^2(D)$,
		\[
		\|M\beta_n-M\beta\|_{L^2(D)}
		\le \|\nabla U_0\|_{L^2(D)}
		\|E(\beta_n-\beta)\|_{L^\infty(D)}\longrightarrow0.
		\]
		The continuous embedding $X_{\rm br}\hookrightarrow L^2(D)$ therefore gives $Z=M\beta$. The closed graph theorem, applied to the underlying real Banach spaces, yields $\|M\beta\|_{X_{\rm br}}\le C_M\|\beta\|_B$. By Theorem~\ref{thm:pullback}, $\beta\mapsto\dot U[\beta]$ is bounded and real-linear. Consequently,
		\[
		U'[\beta]=\dot U[\beta]-M\beta,\quad
		\|U'[\beta]\|_{X_{\rm br}}\le C\|\beta\|_B.
		\]
		This proves the asserted real-linearity and continuity without adding a common-neighborhood regularity assumption.
	\end{proof}

	\section{Statistics of random solutions}\label{section:Statistics}
	
	\subsection{First-order shape Taylor expansion}

	Let $X_K=H^1(K;\mathbb C^d)\times H^1(K;\mathbb C)$. By \eqref{eq:fixed_observation_identity} and Theorem~\ref{thm:pullback}, $U_\epsilon|_K$ is differentiable in $X_K$, with
	\[
	\left.\frac{d}{d\epsilon}U_\epsilon|_K\right|_{\epsilon=0}
	=\dot U[\beta]|_K=:U'[\beta]|_K.
	\]
	This identification requires no additional interface-trace regularity. Under the assumptions of Lemma~\ref{lemma:shape_deriv}, the derivative on $K$ coincides with the restriction of the piecewise Eulerian shape derivative established there.

	\begin{theorem}\label{lemma:taylor}
		Under the assumptions of Theorem~\ref{thm:pullback}, for all sufficiently small $\epsilon>0$ and almost every $\omega\in\Xi$, 
		\begin{equation}\label{eq:taylor_exp}
			U_\epsilon|_K=U_0|_K+\epsilon U'[\beta]|_K+R_\epsilon,
			\quad
			\|R_\epsilon\|_{X_K}\le C\epsilon^2\|\beta\|_B^2,
		\end{equation}
		where $C$ is independent of $\omega$.
	\end{theorem}
	
	\begin{proof}
		Since the perturbation extension vanishes near $\overline K$, \eqref{eq:fixed_observation_identity} gives
		$\widetilde U_\epsilon|_K=U_\epsilon|_K$ and
		$\dot U[\beta]|_K=U'[\beta]|_K$.
		The result follows by restricting \eqref{eq:material_taylor} to $K$.
	\end{proof}

	On the full reference domain, \eqref{eq:material_taylor} gives a first-order expansion of the pulled-back state with linear term $\epsilon\dot U[\beta]$. Under the local regularity assumptions of Lemma~\ref{lemma:shape_deriv}, \eqref{eq:shape_material_relation} yields the following expansion in $L^2(D;\mathbb C^{d+1})$:
	\[
	\widetilde U_\epsilon-U_0-\epsilon U'[\beta]
	=\epsilon(\nabla U_0)V_\beta
	+O_{L^2(D)}\!\left(\epsilon^2\|\beta\|_B^2\right).
	\]
	The transport term $\epsilon(\nabla U_0)V_\beta$ can thus be of first order for an individual realization. Since $V_\beta=E\beta$ with $E$ bounded and real-linear, centering gives $\mathbb E[V_\beta]=E(\mathbb E\beta)=0$. Hence this first-order term vanishes on taking the mean. On $K$ it vanishes for each admissible realization because $V_\beta|_K=0$.

	Throughout this section, $0<\epsilon\le\epsilon_0$, and all statistical errors are measured on the fixed physical observation set $K$. Let $Y_\epsilon=U_\epsilon|_K$, $Y_0=U_0|_K$, and $Z=U'[\beta]|_K=(\mathbf z,z_p)$. By Theorem~\ref{lemma:taylor},
	\begin{equation}\label{eq:statistical_expansion}
		\begin{gathered}
			Y_\epsilon=Y_0+\epsilon Z+R_\epsilon,\quad \mathbb EZ=0,\\
			\|Z\|_{L^2(\Xi;X_K)}\le C,\quad
			\|R_\epsilon\|_{L^2(\Xi;X_K)}\le C\epsilon^2.
		\end{gathered}
	\end{equation}
	For fixed $\epsilon$, the pulled-back forms and loads depend continuously on $\beta$ on the admissible ball. By uniform invertibility, $\beta\mapsto\widetilde U_\epsilon$ is continuous there. Thus $Y_\epsilon$ and the remainder in \eqref{eq:statistical_expansion} are strongly measurable.
	The identity $\mathbb EZ=0$ follows from $\mathbb E\beta=0$ and the fact that the derivative map is bounded and real-linear and therefore commutes with Bochner expectation. The constants may depend on $M_\beta$. Analogous results hold for whole-domain material statistics, with $X$ in place of $X_K$ and with $Y_\epsilon=\widetilde U_\epsilon$, $Y_0=U_0$, and $Z=\dot U[\beta]$.

	\subsection{Mean field approximation}
	
	\begin{theorem}\label{thm:mean}
		Under the assumptions of Theorem~\ref{lemma:taylor}, with $\E[\beta]=0$,

		for all sufficiently small $\epsilon>0$,

		\begin{equation}
			\|\E[\bfu_\epsilon]-\bfu\|_{H^1(K)^d}
			+
			\|\E[p_\epsilon]-p\|_{H^1(K)}
			\le C\epsilon^2 .
		\end{equation}
	\end{theorem}

	\begin{proof}
		Define $T_K:B\to X_K$ by $T_K\beta=\dot U[\beta]|_K$. By Theorem~\ref{thm:pullback}, the map $\beta\mapsto\dot U[\beta]$ from $B$ to $X$ is bounded and real-linear.  Restriction to $K$ is bounded from $X$ to $X_K$, so $T_K$ is bounded and real-linear as well. By \eqref{eq:fixed_observation_identity}, $Z=U'[\beta]|_K=T_K\beta$. Strong measurability of $\beta$ and \eqref{eq:beta_assumption} give Bochner integrability of $Z$ and
		\[
		\mathbb E Z=T_K(\mathbb E\beta)=0 \quad\text{in }X_K.
		\]
		Taking expectations in \eqref{eq:statistical_expansion} now gives $\mathbb E Y_\epsilon-Y_0=\mathbb E R_\epsilon$ in $X_K$. The reverse triangle inequality gives $|\|a\|_{X_K}-\|b\|_{X_K}|\le\|a-b\|_{X_K}$, so the norm is continuous. The triangle inequality gives convexity. 
		{ Applying Jensen's inequality together with the remainder estimate} \eqref{eq:taylor_exp} { yields
			\[
			\|\E R_\epsilon\|_{X_K}
			\le \E\|R_\epsilon\|_{X_K}
			\le C\epsilon^2 \E\|\beta\|_{C^{2,1}(\Gamma)}^2
			\le C\epsilon^2 M_\beta^2.
			\]
			The constant $M_\beta$ can be absorbed into $C$, which completes the proof.}
	\end{proof}
	
	Here, second-order mean accuracy means that the nominal solution alone approximates the mean with an $O(\epsilon^2)$ error.
	
	\subsection{Covariance and variance approximation}
	
	For a complex column vector $\mathbf{a}$, let ${\mathbf{a}}^*=\overline {\mathbf{a}}^{\,T}$, and define the centered field $Y_\epsilon^c=Y_\epsilon-\mathbb EY_\epsilon$. The Hermitian covariance kernel is given by
	\begin{equation}\label{eq:covariance_blocks}
		\mathbf C_{U,\epsilon}(\mathbf x,\mathbf y)
		=\mathbb E\bigl[Y_\epsilon^c(\mathbf x)Y_\epsilon^c(\mathbf y)^*\bigr]
		=\begin{pmatrix}
			\mathbf C_{uu,\epsilon}&\mathbf C_{up,\epsilon}\\
			\mathbf C_{pu,\epsilon}&C_{pp,\epsilon}
		\end{pmatrix}(\mathbf x,\mathbf y).
	\end{equation}
	{ Here $\mathbf C_{uu,\epsilon}=\mathbf C_{u,\epsilon}$ is the displacement covariance, and $C_{pp,\epsilon}=C_{p,\epsilon}$ is the pressure covariance.} 
	In particular,
	\[
	\mathbf C_{up,\epsilon}(\mathbf x,\mathbf y)
	=\mathbb E[\mathbf u_\epsilon^c(\mathbf x)
	\overline{p_\epsilon^c(\mathbf y)}],
	\quad
	\mathbf C_{pu,\epsilon}(\mathbf x,\mathbf y)
	=\mathbf C_{up,\epsilon}(\mathbf y,\mathbf x)^*.
	\]
	These statistics describe complex harmonic amplitudes. The variance of a real time-dependent response can also depend on the chosen phase. For any $w\in L^2(K;\mathbb C^{d+1})$, the kernel satisfies
	\[
	\int_{K\times K}w(\mathbf x)^*\mathbf C_{U,\epsilon}(\mathbf x,\mathbf y)w(\mathbf y)\,d\mathbf x\,d\mathbf y
	=\mathbb E\left|\int_K w(\mathbf x)^*Y_\epsilon^c(\mathbf x)\,d\mathbf x\right|^2\ge0 .
	\]

	Define the scalar displacement and pressure variances by
	\begin{equation}\label{eq:variance_direct}
		V_{u,\epsilon}(\mathbf x)=\mathbb E|\mathbf u_\epsilon^c(\mathbf x)|^2,
		\quad
		V_{p,\epsilon}(\mathbf x)=\mathbb E|p_\epsilon^c(\mathbf x)|^2.
	\end{equation}
	These quantities are well defined almost everywhere in $K$, with the displacement variance equal to the sum of the componentwise variances. This direct definition avoids a diagonal restriction of the covariance kernel. Such a restriction is undefined for a general $L^2(K\times K)$ equivalence class.

	{The covariance kernel of the centered, real-valued interface amplitude is
		\begin{equation}\label{eq:cor_beta}
			\Cor_\beta(\mathbf s,\mathbf t)
			=\operatorname{Cov}(\beta(\mathbf s),\beta(\mathbf t))
			=\mathbb E[\beta(\mathbf s)\beta(\mathbf t)],
			\quad \mathbf s,\mathbf t\in\Gamma.
	\end{equation}}
	The leading response statistics are
	\begin{equation}\label{eq:leading_statistics}
		\mathbf C_Z(\mathbf x,\mathbf y)=\mathbb E[Z(\mathbf x)Z(\mathbf y)^*],
		\quad V_{u,Z}=\mathbb E|\mathbf z|^2,\quad
		V_{p,Z}=\mathbb E|z_p|^2.
	\end{equation}

	\begin{theorem}{ \label{thm:variance}
			Under the assumptions of Theorem~\ref{lemma:taylor},}
		\begin{align}
			\|\mathbf C_{U,\epsilon}-\epsilon^2\mathbf C_Z\|_{L^2(K\times K)}
			&\le C\epsilon^3,\label{eq:covariance_continuum_error}\\
			\|V_{u,\epsilon}-\epsilon^2V_{u,Z}\|_{L^2(K)}
			+\|V_{p,\epsilon}-\epsilon^2V_{p,Z}\|_{L^2(K)}
			&\le C\epsilon^3.\label{eq:variance_continuum_error}
		\end{align}
		The matrix norm in the first estimate is the Frobenius norm, which includes both mixed covariance blocks.
	\end{theorem}
	\begin{proof}
		Let $R_\epsilon^c=R_\epsilon-\mathbb ER_\epsilon$. Then $Y_\epsilon^c=\epsilon Z+R_\epsilon^c$ and $\|R_\epsilon^c\|_{L^2(\Xi;X_K)}\le C\epsilon^2$. For $a,b\in L^2(K;\mathbb C^{d+1})$,
		\[
		\|a(\mathbf x)b(\mathbf y)^*\|_{L^2(K\times K)}
		=\|a\|_{L^2(K)}\|b\|_{L^2(K)}.
		\]
		Expansion of the covariance gives two cross terms, each bounded by
		\[
		\epsilon\|Z\|_{L^2(\Xi;L^2(K))}
		\|R_\epsilon^c\|_{L^2(\Xi;L^2(K))}
		\le C\epsilon^3,
		\]
		and a remainder covariance bounded by $C\epsilon^4$.
		
		For the variance estimate, the continuous embedding $H^1(K)\hookrightarrow L^4(K)$ for $d\le3$ gives
		\[
		\|Z\|_{L^2(\Xi;L^4(K))}\le C,\quad
		\|R_\epsilon^c\|_{L^2(\Xi;L^4(K))}\le C\epsilon^2.
		\]
		For either response block,
		\[
		\mathbb E|\epsilon Z+R_\epsilon^c|^2-\epsilon^2\mathbb E|Z|^2
		=2\epsilon\operatorname{Re}\mathbb E
		(Z\cdot\overline{R_\epsilon^c})
		+\mathbb E|R_\epsilon^c|^2.
		\]
		H\"older's inequality in space and the Cauchy--Schwarz inequality in probability bound the $L^2(K)$ norm of this difference by
		\[
		2\epsilon\|Z\|_{L^2(\Xi;L^4(K))}
		\|R_\epsilon^c\|_{L^2(\Xi;L^4(K))}
		+\|R_\epsilon^c\|_{L^2(\Xi;L^4(K))}^2
		\le C\epsilon^3.
		\]
		The variance estimate is thus independent of a diagonal trace theorem for covariance kernels.
	\end{proof}
	The $O(\epsilon^3)$ estimates above are absolute errors. If $\|\mathbf C_Z\|_{L^2(K\times K)}>0$, the covariance estimate gives an $O(\epsilon)$ relative error with respect to $\epsilon^2\mathbf C_Z$, and also with respect to $\mathbf C_{U,\epsilon}$ for sufficiently small $\epsilon$.

	Earlier studies used first-order second-moment analysis to estimate errors in the mean and variance for Maxwell interface problems~\cite{Hao2018} and elliptic interface problems~\cite{HarbrechtLi2013}. The present results extend this framework to the coupled physical interface problem, linking its shape derivative and numerical approximation to consistent complex-valued statistics.

	\begin{remark}\label{rem:symmetric_inputs}
		Suppose the law of $\beta$ is invariant under $\beta\mapsto-\beta$ and the third derivative of $Y_\epsilon(\beta)$ with respect to $\epsilon$ is uniformly bounded. Writing $Z_2=\partial_\epsilon^2Y_\epsilon(\beta)|_{\epsilon=0}$, we have
		\[
		Y_\epsilon=Y_0+\epsilon Z+\tfrac12\epsilon^2Z_2
		+O_{L^2(\Xi;X_K)}(\epsilon^3).
		\]
		Since $F_{\epsilon,-\beta}=F_{-\epsilon,\beta}$, uniqueness implies that $Z$ is odd and $Z_2$ is even in $\beta$. The order-$\epsilon^3$ cross terms involving $Z$ and $Z_2-\mathbb E Z_2$ therefore vanish. Estimating the remaining terms as in the proof of Theorem~\ref{thm:variance} gives $O(\epsilon^4)$ errors for both the covariance and the variance. The required third-order regularity follows from the uniformly invertible fixed-space equation if its coefficients and sources have uniformly bounded derivatives through order three; the $C^2$ assumption alone is insufficient.
	\end{remark}

	\section{Numerical implementation via low-rank approximation}\label{section:Numerical_Implementation}
	
	The preceding analysis reduces the leading covariance to the covariance of the linear shape derivative. Since $\beta \mapsto (\mathrm{d}\mathbf{u}, \mathrm{d}p)$ is linear, a low-rank representation of ${ \Cor_\beta}$ immediately yields a low-rank representation of the leading covariance.

	The numerical tests use a fixed finite-dimensional random geometry
	\begin{equation}\label{eq:finite_random_geometry}
		\beta(x,\omega)=\sum_{a=1}^q\xi_a(\omega)\psi_a(x),\quad
		\mathbb E\boldsymbol\xi=0,\quad
		Q=\mathbb E[\boldsymbol\xi\boldsymbol\xi^T]\in\mathbb R^{q\times q}.
	\end{equation}
	The smooth modes $\psi_a\in B$, their velocity extensions, and the number of modes $q$ are held fixed during mesh refinement. The random coefficients $\xi_a$ are bounded almost surely.
	Thus
	\begin{equation}\label{eq:finite_mode_covariance}
		\Cor_\beta(\mathbf s,\mathbf t)
		=\sum_{a,b=1}^q\psi_a(\mathbf s)Q_{ab}\psi_b(\mathbf t).
	\end{equation}

	To construct the deterministic interface modes, we apply pivoted Cholesky decomposition to the { coefficient covariance matrix $Q$}. {  A Cholesky column $\boldsymbol\ell_j$ gives the smooth interface mode $\beta_j=\sum_a(\boldsymbol\ell_j)_a\psi_a$.} The algorithm adaptively selects the rank $m$ so that the trace residual of the chosen covariance approximation is at most the prescribed tolerance $\eta_{\rm tol}$. { The residual is measured in the specified Euclidean coefficient space:}

	\begin{equation}\label{eq:coefficient_cholesky_residual}
		R_m=Q-L_mL_m^T\succeq0,\quad\eta_m=\operatorname{tr}R_m,
		\quad L_m=[\boldsymbol\ell_1,\ldots,\boldsymbol\ell_m].
	\end{equation}
	Substituting $Q\approx L_mL_m^T$ into \eqref{eq:finite_mode_covariance} and using the reconstructed modes $\beta_j=\sum_a(\boldsymbol\ell_j)_a\psi_a$ gives
	\begin{equation}\label{eq:chol}
		\Cor_\beta(\bfx,\bfy)
		\approx \sum_{j=1}^{m}\beta_j(\bfx)\beta_j(\bfy),
		\quad \bfx,\bfy\in\Gamma.
	\end{equation}
	In exact arithmetic, $R_m$ in \eqref{eq:coefficient_cholesky_residual} is positive semidefinite. The scalar $\eta_m$ is its matrix trace in coefficient space. The reconstructed modes satisfy
	\[
	\sum_{j=1}^m\|\beta_j\|_{C^{2,1}(\Gamma)}^2
	\le C_\psi\sum_{j=1}^m|\boldsymbol\ell_j|^2
	\le C_\psi\operatorname{tr}Q,
	\]
	where $C_\psi$ is independent of $h,m$ for this fixed basis. We do not assert a uniform bound as \(q\to\infty\).

	\begin{algorithm}[H]

		\caption{Low-rank approximation via pivoted Cholesky decomposition}
		\label{alg:Low_rank}
		\raggedright
		
		\textbf{Input}: { Coefficient} covariance matrix \({ Q\in\mathbb R^{q\times q}}\) and tolerance \(\eta_{\rm tol}>0\).
		
		\textbf{Output}: Low-rank approximation
		\({ Q_m}=\sum_{{ j}=1}^{m}{ \boldsymbol\ell_j}{ \boldsymbol\ell_j}^T\)
		such that
		\(\operatorname{trace}({ Q}-{ Q_m})\le \eta_{\rm tol}\).
		{ Here $\ell_{j,i}=(\boldsymbol\ell_j)_i$.}
		
		Set \(m=0\), \({\bf b}=\operatorname{diag}({ Q})\),
		\(e=\|{\bf b}\|_1\), and
		\({\bf p}=(1,2,\ldots,{ q})\).
		
		\textbf{while} \(e>\eta_{\rm tol}\) \textbf{do}
		
		\hspace{5mm} {  If the largest remaining diagonal entry is zero, stop. }Set \(m=m+1\).
		
		\hspace{5mm} Choose
		\({ j_*\in\operatorname*{arg\,max}_{m\le j\le q}}{\bf b}_{{\bf p}_j}\),
		and swap \({\bf p}_m\) and \({\bf p}_{{ j_*}}\).
		
		\hspace{5mm} { Initialize the new column to zero. }Set
		\({ \ell}_{m,{\bf p}_m}=\sqrt{{\bf b}_{{\bf p}_m}}\).
		
		\hspace{5mm} \textbf{for} \(i=m+1,\ldots,{ q}\) \textbf{do}
		
		\hspace{10mm} Compute
		\[
		{ \ell}_{m,{\bf p}_i}
		=
		\frac{
			{ Q}_{{\bf p}_i,{\bf p}_m}
			-
			\sum_{{ j}=1}^{m-1}
			{ \ell}_{{ j},{\bf p}_i}
			{ \ell}_{{ j},{\bf p}_m}
		}
		{{ \ell}_{m,{\bf p}_m}}.
		\]
		
		\hspace{10mm} Update
		\({\bf b}_{{\bf p}_i}
		=
		\max\{{\bf b}_{{\bf p}_i}-{ \ell}_{m,{\bf p}_i}^2,0\}\),
		where the maximum prevents negative diagonal residuals caused by roundoff.
		
		\hspace{5mm} \textbf{end for}
		
		\hspace{5mm} Compute
		\(e=\sum_{i=m+1}^{{ q}}{\bf b}_{{\bf p}_i}\).
		
		\textbf{end while}
		
	\end{algorithm}

	The output of Algorithm~\ref{alg:Low_rank} is a set of deterministic vectors { \(\{\boldsymbol\ell_j\}_{j=1}^m\), which are the coefficient vectors of the smooth interface modes \(\{\beta_j\}_{j=1}^m\)} in \eqref{eq:chol}. These modes form a low-rank approximation of the interface covariance and are used as inputs for the deterministic shape-derivative problems.

	For each mode $\beta_j$, solve the deterministic { material} derivative problem { \eqref{eq:material_derivative} on the reference domain and restrict the response to $K$}:
	\[
	(\dd\bfu_j,\dd p_j):={ U'[\beta_j]|_K=\dot U[\beta_j]|_K}.
	\]
	Then
	\begin{equation}
		\mathbf C_{\dd\bfu}(\bfx,\bfy)\approx\sum_{j=1}^{m}\dd\bfu_j(\bfx){{{\dd\bfu_j(\bfy)}^{*}}},\quad
		C_{\dd p}(\bfx,\bfy)\approx\sum_{j=1}^{m}\dd p_j(\bfx){\overline{{\dd p_j(\bfy)}}},
		\label{eq:exact_factor_covariance}
	\end{equation}
	and the same solutions give the following mixed covariance blocks
	\begin{equation}\label{eq:mixed_lowrank}
		\mathbf C_{\dd\bfu,\dd p}(\bfx,\bfy)\approx\sum_{j=1}^{m}\dd\bfu_j(\bfx){\overline{{\dd p_j(\bfy)}}},\quad
		\mathbf C_{\dd p,\dd\bfu}(\bfx,\bfy)\approx\sum_{j=1}^{m}\dd p_j(\bfx){\dd\bfu_j(\bfy)^{*}}.
	\end{equation}
	{
		For the finite-mode representation \eqref{eq:finite_random_geometry}, the variance approximations are given almost everywhere on $K$ by
	}
	\begin{equation}\label{eq:variance_du_lowrank}
		{V_{u,\epsilon}}(\bfx)\approx \epsilon^2 \operatorname{tr}\mathbf C_{\dd\bfu}(\bfx,\bfx) \approx \epsilon^2\sum_{j=1}^{m}|\dd\bfu_j(\bfx)|^2, \quad
		{V_{p,\epsilon}}(\bfx)\approx \epsilon^2 C_{\dd p}(\bfx,\bfx) \approx\epsilon^2\sum_{j=1}^{m}|\dd p_j(\bfx)|^2.
	\end{equation}
	Thus, once the nominal state problem has been solved, only $m$ deterministic
	shape-derivative solves are required. The complete low-rank moment approximation is
	summarized in Algorithm~\ref{alg:Low_rank_variance}.

	\begin{algorithm}[htbp]

		\caption{Low-rank moment approximation for the variance of random solutions}
		\label{alg:Low_rank_variance}
		
		\raggedright
		
		\textbf{Input}: Covariance matrix
		{ \(Q\in\mathbb R^{q\times q}\)},
		tolerance \(\eta_{\rm tol}>0\), perturbation amplitude \(\epsilon\), and the
		nominal { discrete} solution \((\bfu_{ h},p_{ h})\).
		
		\textbf{Output}: Leading-order moment approximations of
		\({V_{u,\epsilon}}({\bf x})\) and \({V_{p,\epsilon}}({\bf x})\).
		
		\medskip
		
		Step 1. \textbf{Low-rank approximation of the interface covariance.}\\
		\hspace{13mm} Compute
		\(\Cor_\beta\approx \sum_{j=1}^m \beta_j\otimes\beta_j\)
		by Algorithm~\ref{alg:Low_rank}{  applied to $Q$, with $\beta_j=\sum_a\psi_a(\boldsymbol\ell_j)_a$}.
		
		\medskip
		
		Step 2. \textbf{Solve the shape derivative problems.}\\
		\hspace{13mm} For each mode \(\beta_j\), compute
		\((\dd\bfu_{j{ ,h}},\dd p_{j{ ,h}}){ =\dot U_h[V_{\beta_j,h}]|_K}\)
		by solving {  the the differentiated discrete weak problem \eqref{eq:discrete_sensitivity} with $V_{\beta_j,h}=I_hV_{\beta_j}$. On $K$, the resulting discrete material sensitivity agrees with the Eulerian sensitivity when $V_{\beta_j,h}|_K=0$.} 
		
		\medskip
		
		Step 3. \textbf{Construct the leading covariance of the perturbations.}\\
		\hspace{13mm} Compute
		\((\mathbf C_{\dd\bfu{ ;h,m}}({\bf x},{\bf y}),C_{\dd p{ ;h,m}}({\bf x},{\bf y}))\)
		using \eqref{eq:exact_factor_covariance} {  with the discrete responses $(\dd\bfu_{j,h},\dd p_{j,h})$}.
		
		\medskip
		
		Step 4. \textbf{Reconstruct the leading-order variance.}\\
		\hspace{13mm} Compute
		\((\epsilon^2\sum_{j=1}^{m}|\dd\bfu_{j{ ,h}}({\bf x})|^2,
		\epsilon^2\sum_{j=1}^{m}|\dd p_{j{ ,h}}({\bf x})|^2)\),
		which gives the leading-order approximation of
		\(({V_{u,\epsilon}}({\bf x}),{V_{p,\epsilon}}({\bf x}))\)
		according to \eqref{eq:variance_du_lowrank}.
		
	\end{algorithm}

	Algorithm~\ref{alg:Low_rank_variance} requires one nominal Biot solve and $m$ deterministic sensitivity solves, where $m$ is the retained covariance rank. The sensitivity solves reuse the nominal factorization on the fixed reference configuration. The total cost additionally depends on covariance factorization, mode construction, sensitivity right-hand-side assembly, observation evaluation, and the chosen representation of the output covariance. A direct Monte Carlo estimator requires one random-domain solve per sample; the relative cost therefore depends on the target accuracy and on these implementation costs.

	The next lemma quantifies the propagation of the Cholesky truncation error through the linear shape-derivative map.
	
	\begin{lemma}\label{lem:low_rank_consistency}\label{lem:covtrace}
		{ Let $\mathcal S:\mathbb C^q\to X_K$ be the complex-linear extension of the map sending a real coefficient vector $\boldsymbol c$ to $U^\prime[\sum_a c_a\psi_a]|_K$. Set $Q_m=L_mL_m^T$. With $\eta_m=\operatorname{tr}(Q-Q_m)$ as in \eqref{eq:coefficient_cholesky_residual},} we have
		\begin{equation}\label{eq:covtrace_bound}
			\|\mathcal S { Q} \mathcal S^*-\mathcal S { Q_m}\mathcal S^*\|_{\mathcal L_1(X_K)}
			\le \|\mathcal S\|^2\eta_m.
		\end{equation}
		Consequently, the low-rank contribution to the leading covariance error is $O(\epsilon^2\eta_m)$. { The same order holds for the pressure variance and displacement trace variance in $L^2(K)$, for $d=2,3$.}
	\end{lemma}
	\begin{proof}
		{  With the product $H^1$ inner product on $X_K$, let $\mathcal S^*$ denote the Hilbert adjoint. The derivative covariance operator on $X_K$ is}
		\[
		{ \mathcal C_Z}=\mathcal S { Q}\mathcal S^*.
		\]
		The rank-$m$ approximation gives ${ \mathcal C_{Z,m}}=\mathcal S { Q_m}\mathcal S^*$.  Hence
		\[
		{ \mathcal C_Z}-{ \mathcal C_{Z,m}}=\mathcal S({ Q}-{ Q_m})\mathcal S^*.
		\]
		The trace-class ideal property gives
		\[
		\|\mathcal S({ Q}-{ Q_m})\mathcal S^*\|_{\mathcal L_1(X_K)}
		\le \|\mathcal S\|_{\mathcal L({ \mathbb C^q},X_K)}^2
		\|{ Q}-{ Q_m}\|_{\mathcal L_1({ \mathbb C^q})}.
		\]
		{  Since $R_m$ is positive semidefinite, its trace norm equals $\operatorname{tr}R_m=\eta_m$. This proves \eqref{eq:covtrace_bound}. Let $\iota:X_K\hookrightarrow L^2(K;\mathbb C^{d+1})$ be the natural embedding, so that the physical covariance kernel represents $\iota\mathcal C_Z\iota^*$. Decompose $R_m=\sum_a r_a\mathbf e_a\mathbf e_a^*$, where $r_a\ge0$ and $\{\mathbf e_a\}$ is orthonormal. Then the Frobenius norm of the residual kernel satisfies
			\[
			\left\|\sum_a r_a(\mathcal S\mathbf e_a)(\mathbf x)
			(\mathcal S\mathbf e_a)(\mathbf y)^*\right\|_{L^2(K\times K)}
			\le\sum_a r_a\|\mathcal S\mathbf e_a\|_{L^2(K)}^2
			\le C\|\mathcal S\|^2\eta_m.
			\]
			To bound the variance directly, use the same decomposition. For the pressure or displacement projection $\pi$, Sobolev embedding gives
			\[
			\left\|\sum_a r_a|\pi\mathcal S\mathbf e_a|^2\right\|_{L^2(K)}
			\le\sum_a r_a\|\pi\mathcal S\mathbf e_a\|_{L^4(K)}^2
			\le C\|\mathcal S\|^2\eta_m.
			\]
			Here $H^1(K)\hookrightarrow L^4(K)$ for $d=2,3$. Multiplication by $\epsilon^2$ proves the variance bound without a kernel diagonal restriction.}
	\end{proof}

	\section{Error estimates}\label{section:Error_Estimates}
	
	The error has three distinct sources: the shape Taylor truncation, the finite element discretization, and the low-rank truncation of the perturbation covariance.  Separating these terms is important because the low-rank tolerance is independent of the interface perturbation amplitude $\epsilon$. We first describe the finite element setting used for the nominal state and the shape-derivative problems, and then state the fully discrete estimate.
	
	Let $\mathcal T_h$ be an interface-fitted, shape-regular triangulation {  or tetrahedralization}.  The deterministic state is approximated by conforming finite elements.

	We use continuous $\mathbb P_2$ elements for displacement and continuous $\mathbb P_1$ elements for pressure. We then differentiate the discrete weak equations. This gives the discrete material derivative. It agrees with the discrete Eulerian derivative on the fixed observation set $K$ when the interpolated velocity vanishes on $K$. No simultaneous lifting of the four interface jumps is used.
	For a mesh that resolves the interface exactly, define the complex finite element space
	\[
	X_h=\left\{(\mathbf v_h,q_h)\in X:\;
	\mathbf v_h|_T\in\mathbb P_2(T;\mathbb C)^d,\quad
	q_h|_T\in\mathbb P_1(T;\mathbb C),\quad
	T\in\mathcal T_h\right\}.
	\]
	Here $\mathbb P_k(T;\mathbb C)$ denotes the polynomials of total degree at most $k$ with complex coefficients. The nodal basis functions are real-valued, and the nodal coefficients are complex.
	
	The nominal discrete problem is to find $U_h\in X_h$ such that
	\begin{equation}\label{eq:discrete_state}
		b_0(U_h,W_h)=\ell_0(W_h),\quad \forall\,W_h\in X_h.
	\end{equation}
	Test functions are conjugated in every sesquilinear form. The name of the element pair does not provide a separate stability argument; the fixed-parameter justification is given below.

	For a smooth perturbation velocity $V=V_\beta$, let $V_h=I_hV$ denote its continuous piecewise affine nodal interpolant, and write $\dot U[V]=\dot U[\beta]$. 
	For the continuous forms, $\dot b_0[V]$ and $\dot\ell_0[V]$ denote $b'_0[V]$ and $\ell'_0[V]$ in \eqref{eq:form_derivative} and \eqref{eq:load_derivative}, respectively.
	Differentiating the discrete equations pulled back by $\Phi_{t,h}=I+tV_h$, with fixed material labels and reference coefficients, yields
	\begin{equation}\label{eq:discrete_sensitivity}
		b_{0,h}(\dot U_h[V_h],W_h)
		=\dot\ell_{0,h}[V_h](W_h)
		-\dot b_{0,h}[V_h](U_h,W_h),
		\quad W_h\in X_h.
	\end{equation}
	For exact geometry and exact integration, $b_{0,h}=b_0|_{X_h\times X_h}$. On a polygonal approximation, $U_h$ instead solves the nominal equations assembled on that same approximate geometry, and both derivatives in \eqref{eq:discrete_sensitivity} refer to its assembled forms. In either case, the linear system is $A_h\dot{\mathbf U}_h=\dot{\mathbf F}_h-\dot A_h\mathbf U_h$. All sensitivity systems therefore share the nominal matrix $A_h$.

	This differentiates the fitted Biot residual consistently with its discrete geometry. Related work analyzes shape derivatives of linear functionals in structural topology optimization~\cite{KleinNairYano2022} and discrete shape functionals on unfitted meshes~\cite{ShahanWalker2026}. The estimates below concern the fitted Biot state and its material derivative.

	For $G_T=\nabla V_h|_T$, the elementwise geometric derivatives satisfy
	\begin{equation}\label{eq:element_geometry_derivative}
		\left.\frac{d}{dt}\det(I+tG_T)\right|_{t=0}=\operatorname{tr}G_T,
		\quad
		\left.\frac{d}{dt}(I+tG_T)^{-1}\right|_{t=0}=-G_T.
	\end{equation}
	We differentiate the gradient transformations, the volume factors, and the compositions of the fixed physical sources with the deformation map. The last derivatives involve the spatial gradients of the sources, as shown in \eqref{eq:load_derivative}. For the constant sources and piecewise constant coefficients used in the experiments, four-point Duffy--Gauss quadrature in each transformed coordinate integrates the element contributions and their derivatives exactly. The quadrature was verified against the simplex moment formula
	\begin{equation}\label{eq:simplex_moments}
		\int_T\prod_{a=0}^d\lambda_a^{r_a}\,dx
		=d!\,|T|\,\frac{\prod_{a=0}^d r_a!}{(d+\sum_a r_a)!},
	\end{equation}
	where the $\lambda_a$ are barycentric coordinates and the $r_a$ are nonnegative integers. The highest polynomial degree on an affine simplex is four, from the displacement mass term. After the Duffy transform, each coordinate degree is at most $4+d-1\le6$, including the Jacobian. Four-point Gauss integration is exact through degree seven in each coordinate.

	The velocity $V$ vanishes on the fixed observation set $K$ by the choice of extension in Section~\ref{sec:stochastic_geometry}. If the perturbation velocity vanishes at all vertices of elements intersecting $K$, its piecewise affine interpolant $V_h$ also vanishes on $K$. Consequently, the material and Eulerian shape sensitivities coincide on $K$.

	The planar interface in the 3D test is represented exactly. For the polygonal approximation of the 2D circular interface, assume a uniformly bi-Lipschitz geometric lift that is the identity near $K$. Let $\widehat X_h\subset X$ denote the conforming lifted spaces satisfying
	\[
	\inf_{Z_h\in\widehat X_h}\|Z-Z_h\|_X\to0
	\quad\text{for every }Z\in X.
	\]
	Assume that the lifted forms and loads admit bounded extensions to the continuous spaces and satisfy
	\begin{equation}\label{eq:audit_geometry_consistency}
		\begin{aligned}
			\delta_{\mathrm{geo},h}:={}&
			\|b_h^\ell-b_0\|_{\mathcal L(X,X^*)}
			+\|\ell_h^\ell-\ell_0\|_{X^*}\\
			&+\sup_{\substack{V\in\mathcal V\\\|V\|_{W^{2,\infty}}\le1}}
			\left(
			\|\dot b_h^\ell[V]-\dot b_0[V]\|_{\mathcal L(X,X^*)}
			+\|\dot\ell_h^\ell[V]-\dot\ell_0[V]\|_{X^*}
			\right)
			\longrightarrow0.
		\end{aligned}
	\end{equation}
	Here $\mathcal V=\operatorname{span}_{\mathbb R} \{V_{\psi_a}\}_{a=1}^q$. A smooth physical velocity is restricted to the polygonal mesh before interpolation. We assume that the lifted differentiated forms and loads are uniformly Lipschitz continuous in the velocity with respect to the $W^{1,\infty}$ norm, including under replacement of $V$ by $I_hV$. These lifting and consistency properties remain unproved for the polygonal circle meshes, so the subsequent estimates are conditional on them. Quadrature, sampling, and algebraic solver errors are excluded.

	Geometry approximation is also important in unfitted interface methods \cite{ChenLiu2023,HeimannLehrenfeld2025}. For the present fitted scheme, both the nominal operator and its shape derivative require the lifted consistency assumptions in \eqref{eq:audit_geometry_consistency}.

	\begin{theorem}
		\label{thm:discrete}
		Under Assumption~\ref{ass:data} and the coefficient conditions of Section~\ref{section:Model}, assume $\inf_{Z_h\in X_h}\|Z-Z_h\|_X\to0$ for every $Z\in X$.
		For sufficiently small $h$, the discrete state problem \eqref{eq:discrete_state} and the sensitivity problem \eqref{eq:discrete_sensitivity} are uniquely solvable. 
		On exact geometry with exact integration and the same velocity in both differentiated problems,
		\begin{align}
			\|U_0-U_h\|_X&\le C\inf_{Z_h\in X_h}\|U_0-Z_h\|_X,
			\label{eq:state_quasioptimal}\\
			\|\dot U[V]-\dot U_h[V]\|_X
			&\le C\left\{\inf_{Z_h\in X_h}\|\dot U[V]-Z_h\|_X
			+\|V\|_{W^{1,\infty}}\|U_0-U_h\|_X\right\}.
			\label{eq:sensitivity_quasioptimal}
		\end{align}
		Replacing $V$ by $V_h$ in the discrete sensitivity problem adds an error bounded by $C\|V-V_h\|_{W^{1,\infty}} \le Ch\|V\|_{W^{2,\infty}}$. If $U_0$ and the sensitivities for the $q$ fixed basis velocities are piecewise $H^2$, the product-$H^1$ errors are $O(h)$ on exact fitted affine meshes, uniformly for $|\boldsymbol c|\le1$. For approximate geometry satisfying \eqref{eq:audit_geometry_consistency}, the lifted estimates acquire an additional $C\delta_{\mathrm{geo},h}$. The rate $O(h+\delta_{\mathrm{geo},h})$ further requires $O(h)$ best approximation in the lifted spaces. The constants may depend on the fixed physical parameters and modes.
	\end{theorem}

	\begin{proof}
		Lemma~\ref{lem:fredholm} gives the principal-form representation $I+\mathcal K$, with $\mathcal K$ compact. The corresponding orthogonal projections $P_h$ satisfy $P_h\to I$ strongly, hence $\|(I-P_h)\mathcal K\|\to0$. Non-resonance gives eventual uniform invertibility of $I+P_h\mathcal K$ and \eqref{eq:state_quasioptimal}.
		Subtracting the differentiated equations gives
		\[
		b_0(\dot U[V]-\dot U_h[V],W_h)
		=-\dot b_0[V](U_0-U_h,W_h).
		\]
		Using $|\dot b_0[V](Z,W)|\le C\|V\|_{W^{1,\infty}}\|Z\|_X\|W\|_X$, subtract an arbitrary approximant of $\dot U[V]$ and apply discrete stability to obtain \eqref{eq:sensitivity_quasioptimal}. Uniform stability applied to the difference of the two discrete sensitivity equations bounds the resulting change in the discrete material derivative by $C\|V-I_hV\|_{W^{1,\infty}}$. Elementwise interpolation of $U_0$ and of the sensitivities for each fixed basis direction gives the stated rate. Linearity and Cauchy--Schwarz in the $q$ coefficients make the bound uniform on the Euclidean unit ball.
		
		The same argument gives an inf-sup constant $\gamma>0$ on $\widehat X_h$.
		For $\|b_h^\ell-b_0\|<\gamma/2$, the lifted scheme retains stability with constant $\gamma/2$. Insert arbitrary lifted approximants into the state and sensitivity error equations and apply this stability bound. Taking the infimum over the approximants gives the best-approximation terms. The consistency bound \eqref{eq:audit_geometry_consistency} controls the additional form and load residuals by $C\delta_{\mathrm{geo},h}$ on the fixed velocity unit sphere.
		This proves the conditional curved-geometry estimates.
	\end{proof}

	\begin{theorem}\label{thm:6.1}\label{thm:discrete_moments}
		Assume the hypotheses of Theorems~\ref{lemma:taylor}--\ref{thm:variance}. {  Assume the fixed smooth-mode representation \eqref{eq:finite_random_geometry} and the hypotheses of Theorem~\ref{thm:discrete}. Define $\mathcal S_h\boldsymbol c=\dot U_h[I_h\sum_a c_aV_{\psi_a}]|_K$ for real $\boldsymbol c$ and extend $\mathcal S_h$ by complex linearity. Then $(\mathrm d\mathbf u_{j,h},\mathrm dp_{j,h})=\mathcal S_h\boldsymbol\ell_j$.
			Set
			\[
			e_{0,h}=\|U_0|_K-U_h|_K\|_{X_K},\quad
			e_{1,h}=\|\mathcal S-\mathcal S_h\|_{\mathcal L(\mathbb C^q,X_K)}.
			\]
			Define the discrete variance fields
			\[
			V_{u,h,m}:=\epsilon^2\sum_{j=1}^m|\mathrm d\mathbf u_{j,h}|^2,\quad
			V_{p,h,m}:=\epsilon^2\sum_{j=1}^m|\mathrm dp_{j,h}|^2,
			\]
			and let $\eta_m=\operatorname{tr}R_m$ as in \eqref{eq:coefficient_cholesky_residual}.} Then, {  on the fixed observation set $K$},
		\begin{align}
			\|\E[\bfu_\epsilon]-\bfu_h\|_{H^1(K)^d}
			+
			\|\E[p_\epsilon]-p_h\|_{H^1(K)}
			&\le C(\epsilon^2+{ e_{0,h}}), \label{eq:fully_discrete_mean}\\
			\|{V_{u,\epsilon}}-{V_{u,h,m}}\|_{L^2(K)}
			+
			\|{V_{p,\epsilon}}-{V_{p,h,m}}\|_{L^2(K)}
			&\le C(\epsilon^3+\epsilon^2 { e_{1,h}\operatorname{tr}Q}+\epsilon^2\eta_m). \label{eq:fully_discrete_variance}
		\end{align}
		The constant $C$ is independent of $\epsilon$, $h$, and $m$ for sufficiently small $\epsilon$ and sufficiently fine interface-fitted meshes.
		{  Theorem~\ref{thm:discrete} gives $e_{0,h}+e_{1,h}\le Ch$ under piecewise $H^2$ regularity on exact geometry, or $C(h+\delta_{\mathrm{geo},h})$ under its curved-geometry and lifted-approximation conditions.}
	\end{theorem}
	\begin{proof}
		We write $U_{ 0}=(\bfu_{ 0},p_{ 0})$, $U_h=(\bfu_h,p_h)$, and use the product norm on $X_K$.  The mean estimate is obtained by a triangle inequality:
		\[
		\|\E U_\epsilon-U_h\|_{X_K}
		\le
		\|\E U_\epsilon-U_{ 0}\|_{X_K}
		+\|U_{ 0}-U_h\|_{X_K}.
		\]
		The first term is $O(\epsilon^2)$ by Theorem~\ref{thm:mean}; the second is { $e_{0,h}$ by definition}.  This proves \eqref{eq:fully_discrete_mean}.
		
		For the variance, introduce three intermediate quantities.  Let $V_\epsilon$ be the exact variance, $V^{(1)}$ the exact leading variance, $V^{(1)}_m$ the leading variance obtained from the rank-$m$ covariance approximation, and $V^{(1)}_{h,m}$ its finite element approximation.  Then
		\[
		\|V_\epsilon-V^{(1)}_{h,m}\|
		\le
		\|V_\epsilon-V^{(1)}\|+
		\|V^{(1)}-V^{(1)}_m\|+
		\|V^{(1)}_m-V^{(1)}_{h,m}\|.
		\]
		{ Each error term on the right is measured in $L^2(K)$ for either the scalar displacement variance or the pressure variance.} The first term is $O(\epsilon^3)$ by Theorem~\ref{thm:variance}.  The second term is $O(\epsilon^2\eta_m)$ by Lemma~\ref{lem:low_rank_consistency}.
		
		It remains to estimate the finite element part.  For the pressure component, for example,
		\[
		{
			{\bigl|}|\dd p_j|^2-|\dd p_{j,h}|^2{ \bigr|\le}  |\dd p_j-\dd p_{j,h}|( |\dd p_j|+|\dd p_{j,h}|).}
		\]
		The same {inequality} holds component-wise for the displacement trace variance. By H\"older's inequality and { the embedding $H^1(K)\hookrightarrow L^4(K)$ for $d=2,3$},
		\[\begin{aligned}
			\||\dd p_j|^2-|\dd p_{j,h}|^2\|_{L^2(K)}
			&\le \|\dd p_j-\dd p_{j,h}\|_{L^4(K)}
			{
				\bigl(\|\dd p_j\|_{L^4(K)}+\|\dd p_{j,h}\|_{L^4(K)}\bigr)}\\
			&\le C { e_{1,h}(\|\mathcal S\|+\|\mathcal S_h\|)|\boldsymbol\ell_j|^2}.
		\end{aligned}\]
		Summing over the selected Cholesky modes and using { $\sum_j|\boldsymbol\ell_j|^2\le\operatorname{tr}Q$ and the uniform bounds on $\mathcal S,\mathcal S_h$} gives
		\[
		\|V^{(1)}_m-V^{(1)}_{h,m}\|_{L^2(K)}
		\le C\epsilon^2 { e_{1,h}\operatorname{tr}Q}.
		\]
		Combining the three bounds gives \eqref{eq:fully_discrete_variance} for the displacement and pressure variances.
	\end{proof}

	\section{Numerical experiments}\label{section:Numerical_experiments}

	{ This section presents numerical experiments for a circular interface in two dimensions and a planar interface in three dimensions, the perturbations are supported away from the edges of the planar interface. We examine the discrete shape Taylor remainder, spatial approximation, covariance compression, and sampling error. The circular interface is approximated by a fitted polygon, so the spatial errors also include geometric approximation. Fixed-mesh perturbation tests assess discrete Taylor convergence rather than continuum spatial accuracy.}

	\subsection{Experimental setup and accuracy tests}
	\label{sec:exp_convergence}
	
	This subsection presents the benchmark setup and examines the numerical implementation and spatial resolution. The spatial refinement study evaluates the nominal solution and interface sensitivities by comparison with finite element solutions on finer reference meshes.

	The computational domain is $D=(-1,1)^d$, and the dimensionless benchmark parameters are scaled according to \eqref{eq:nondimensionalization} and listed in Table~\ref{tab:parameters}.
	\begin{table}[H]

		\centering
		\caption{Material coefficients used in all numerical experiments. Material 1 occupies the interior of the circle in 2D and the lower half of the computational domain in 3D.}
		\label{tab:parameters}
		\begin{tabular}{lrrrrrr}
			\toprule
			Material & $\lambda$ & $\mu$ & $\alpha$ & $c$ & $\kappa$ & $\rho$\\
			\midrule
			$D_1$ & 4 & 2 & 0.9 & 0.6 & 0.08 & 1.2\\
			$D_2$ & 2 & 1 & 0.7 & 1.0 & 0.30 & 1.0\\
			\bottomrule
		\end{tabular}
	\end{table}
	In both dimensions, we set $\theta=0.8$ and $g=0.6$. We use $\mathbf f=(1,0.25)^T$ in 2D and $\mathbf f=(1,0.25,-0.15)^T$ in 3D. These sources are constant, real, nonzero, and fixed in physical coordinates. Homogeneous Dirichlet conditions are imposed on both displacement and pressure. The solution is complex-valued due to the $i\theta$ terms in the storage and coupling equations. The sufficient non-resonance condition of Lemma~\ref{lem:low_frequency_nonresonance} holds, since $\theta^2\rho_{\max}=0.768<\mu_{\min}\lambda_1(D)$, where $\mu_{\min}\lambda_1(D)=\pi^2/2$ in 2D and $3\pi^2/4$ in 3D. This bound remains valid under admissible perturbations of the internal material interface.

\begin{figure}[htp]

		\centering
		\includegraphics[width=\linewidth]{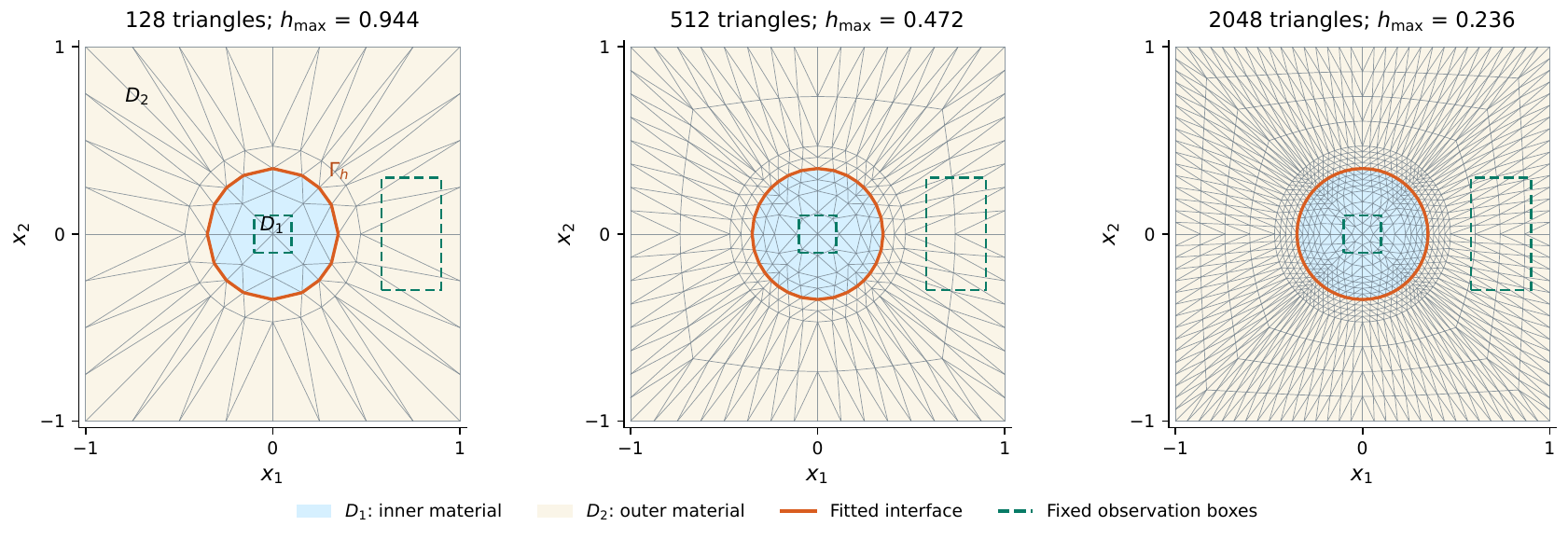}
		\caption{Interface-fitted triangular meshes at three refinement levels.}
		\label{fig:actual_mesh2d}
	\end{figure}

\begin{figure}[H]

		\centering
		\includegraphics[width=0.75\linewidth]{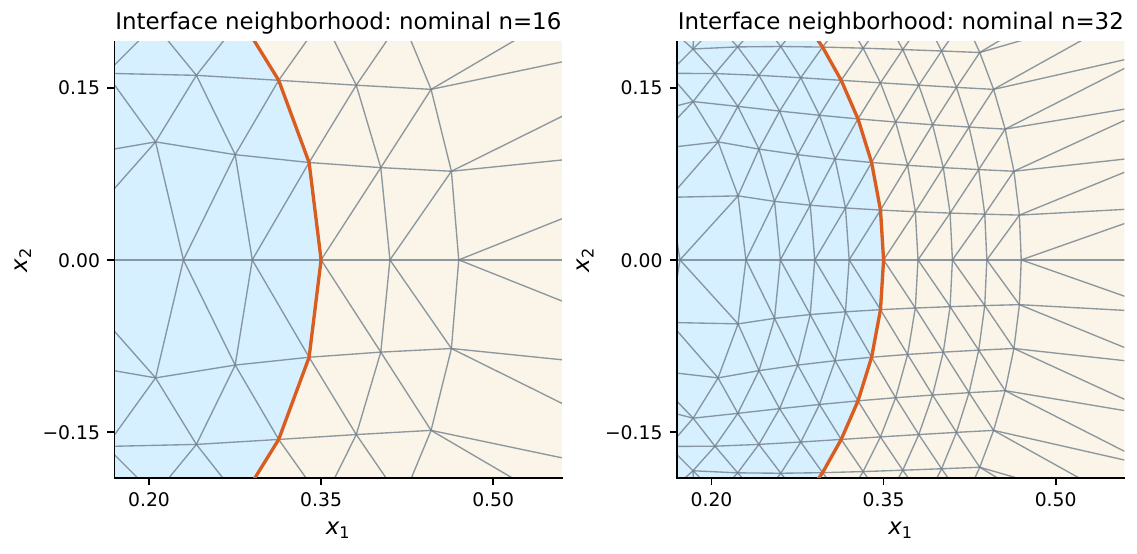}
		\caption{Enlarged views of the interface-fitted triangular meshes at two refinement levels.}
		\label{fig:interface_zoom}
	\end{figure}

	In 2D, $\Gamma$ is the circle $|\mathbf x|=0.35$. The fixed observation set is
	\begin{equation}\label{eq:numerical_K2}
		K_2=(-0.10,0.10)^2
		\ \cup\ (0.58,0.90)\times(-0.30,0.30).
	\end{equation}
	We use interface-fitted triangular meshes with straight edges, resulting in a polygonal approximation of the circular interface. Figure~\ref{fig:actual_mesh2d} shows three levels of mesh refinement, while Figure~\ref{fig:interface_zoom} provides a close-up view of the mesh near the interface. The resulting geometric consistency error is accounted for by $\delta_{\mathrm{geo},h}$.

\begin{figure}[htp]

		\centering
		\includegraphics[width=0.5\linewidth]{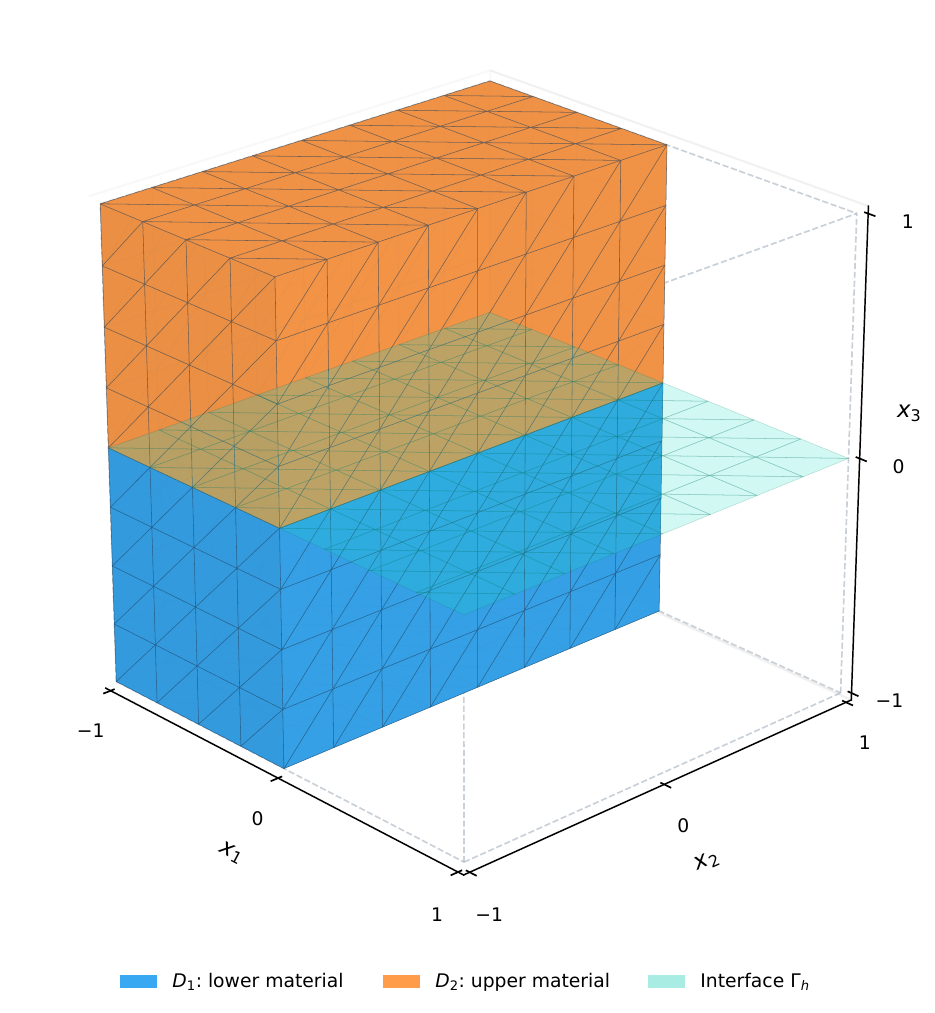}
		\caption{Interface-fitted tetrahedral mesh with 3072 elements, showing the two material subdomains separated by the interface $x_3=0$ in the half-domain $x_1\leq 0$.}
		\label{fig:actual_mesh3d}
	\end{figure}

	In 3D, $\Gamma=\{x_3=0\}\cap D$, $D_1=\{x_3<0\}\cap D$, and
	\begin{equation}\label{eq:numerical_K3}
		K_3=(-0.5,0.5)^2\times
		\bigl((-0.85,-0.5)\cup(0.5,0.85)\bigr).
	\end{equation}
	Each Cartesian cube is subdivided into six Freudenthal tetrahedra. Figure~\ref{fig:actual_mesh3d} shows a half-domain view of the resulting mesh. The nominal planar interface is represented exactly, although its perturbations are generally nonplanar.

	The velocity extensions use the smooth cutoff function
	\[
	\chi(t)=
	\begin{cases}
		\exp[-t^2/(1-t^2)],& |t|<1,\\
		0,& |t|\ge1.
	\end{cases}
	\]
	In 2D, a circular interface perturbation $\beta(\varphi)$ is extended as $\chi((r-0.35)/0.12)\beta(\varphi)\mathbf e_r$, with support confined to an annulus away from the origin. In 3D, define
	\[
	s_*(t)=
	\begin{cases}
		0,&t\le0,\\
		\displaystyle\frac{\exp(-1/t)}
		{\exp(-1/t)+\exp(-1/(1-t))},&0<t<1,\\
		1,&t\ge1,
	\end{cases}
	\quad
	a(t)=s_*((0.9-|t|)/0.1).
	\]
	The 3D surface modes include the factor $a(x_1)a(x_2)$ and are extended in the $\mathbf e_3$ direction using the cutoff $\chi(x_3/0.26)$, ensuring that the velocity vanishes near the interface edges and $\partial D$. Under piecewise affine nodal interpolation, each $\mathbb{P}_2$ edge midpoint receives the average displacement of its endpoints. All vertices of elements intersecting $K$ remain fixed, ensuring comparisons at identical physical observation coordinates.


	All deterministic solvers use an interface-resolved finite element mesh.

	The implementation employs continuous vector-valued ${\mathbb{P}}_2$ elements for displacement, ${\mathbb{P}}_1$ elements for pressure, polynomial-exact element integration, and a general complex sparse LU solver. The four interface transmission conditions are incorporated through the weak formulation and its differentiated discrete system, without pointwise enforcement of conormal jumps. A single nominal matrix factorization is reused for all sensitivity right-hand sides. Since a direct solver is used, no iterative stopping tolerance is required. The maximum relative algebraic residual over the free degrees of freedom (DOFs) is \ensuremath{1.97\times10^{-11}}, and the minimum ratio of perturbed to reference element Jacobian determinants is 0.330. All tested geometric mappings preserve orientation.

	For spatial error comparisons, let $U_{0,h}$ and $U_{0,\mathrm{ref}}$ denote the nominal solutions on the current and reference meshes, respectively, and let $S_{h,j}$ and $S_{\mathrm{ref},j}$ denote the corresponding sensitivities for $j=1,2$. On the observation set $K$, define
	\[
	\|(\mathbf w,r)\|_{s,K}^2
	=\|\mathbf w\|_{H^s(K)^d}^2+\|r\|_{H^s(K)}^2,
	\quad s=0,1,
	\]
	with $H^0(K)=L^2(K)$. The spatial errors are measured by
	\[
	e_U^{H^s}
	=\|U_{0,h}-U_{0,\mathrm{ref}}\|_{s,K},
	\quad
	e_S^{H^s}
	=\left(\sum_{j=1}^2
	\|S_{h,j}-S_{\mathrm{ref},j}\|_{s,K}^2
	\right)^{1/2}.
	\]
	These errors are absolute norm differences from fine-mesh finite element reference solutions.
	
\begin{figure}[H]

		\includegraphics[width=0.94\linewidth]{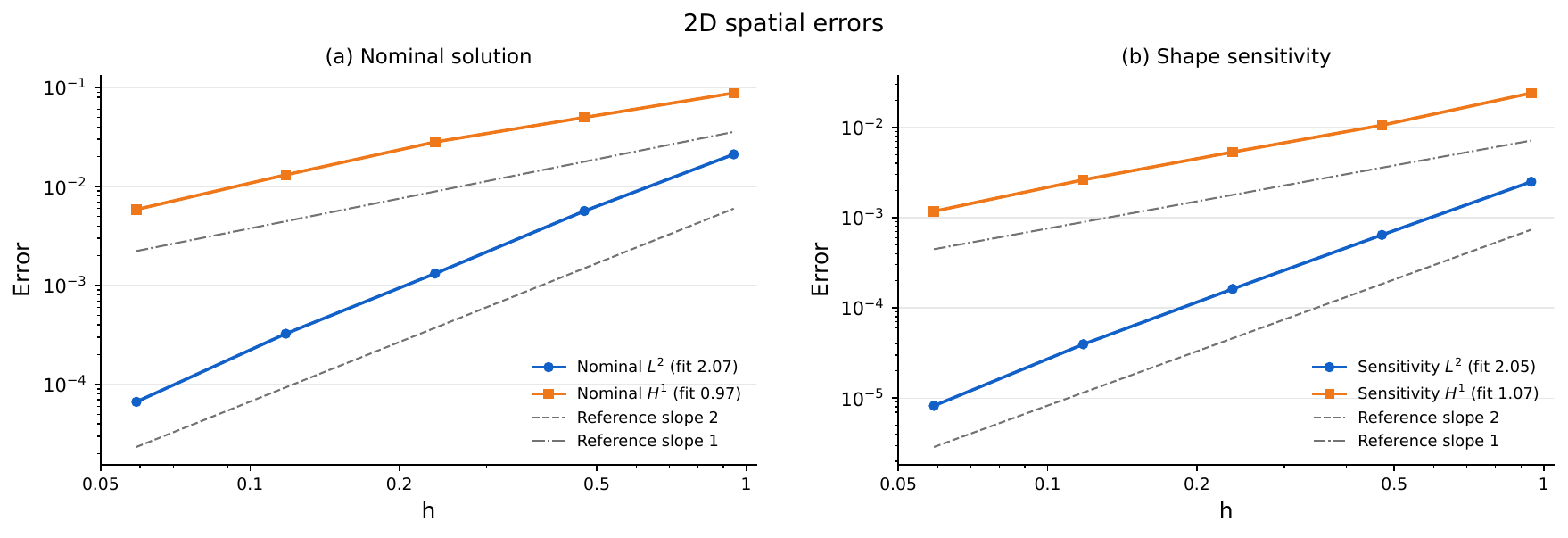}
		\caption{2D spatial errors of the nominal solution (left) and shape sensitivities (right) in the $L^2(K)$ and $H^1(K)$ norms, measured relative to the $n=256$ reference solution.}
		\label{fig:spatial2d}
	\end{figure}

\begin{figure}[H]

		\includegraphics[width=0.94\linewidth]{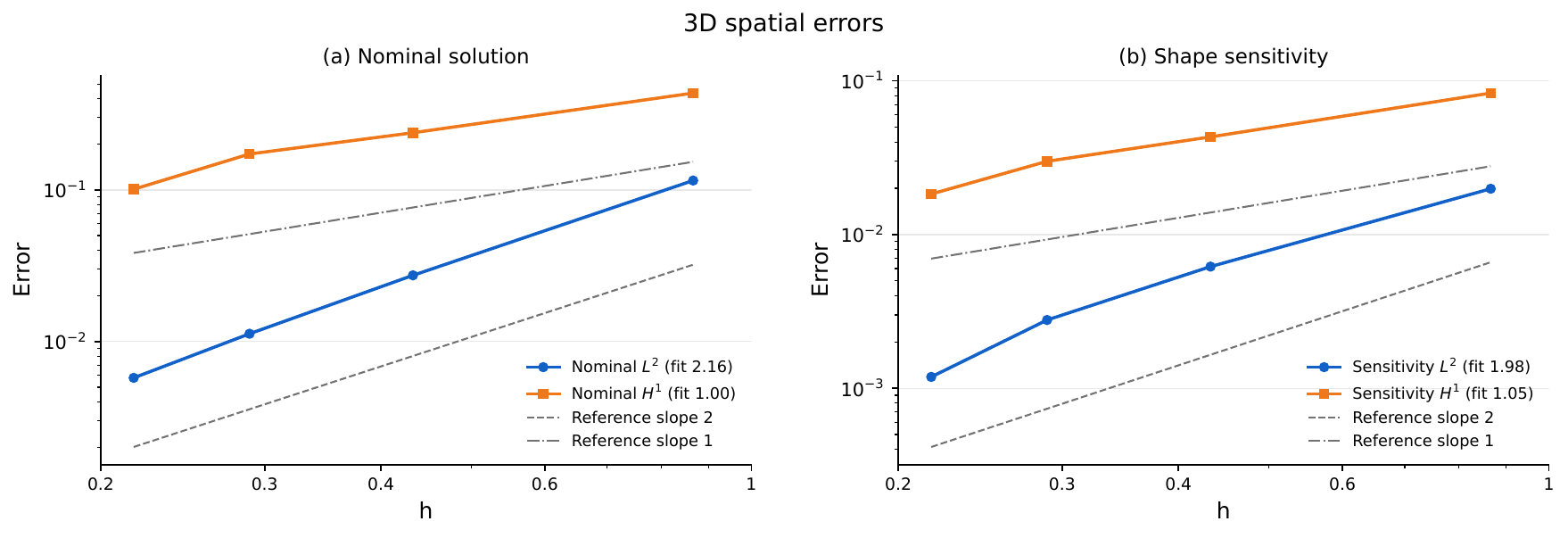}
		\caption{3D spatial errors of the nominal solution (left) and shape sensitivities (right) in the $L^2(K)$ and $H^1(K)$ norms, measured relative to the $n=32$ reference solution.}
		\label{fig:spatial3d}
	\end{figure}

	Table~\ref{tab:convergence} and Figure~\ref{fig:spatial2d} give the 2D results. The 2D reference mesh ($n=256$) has 131072 triangles and 587267 free DOFs. At $n=128$, the nominal and sensitivity $L^2$ errors are $6.669\times10^{-5}$ and $8.218\times10^{-6}$, respectively. Table~\ref{tab:convergence3d} and Figure~\ref{fig:spatial3d} give the 3D results. The 3D reference mesh ($n=32$) has 196608 tetrahedra and 779932 free DOFs. At $n=16$, the corresponding $L^2$ errors are $5.765\times10^{-3}$ and $1.185\times10^{-3}$, respectively.
	In both dimensions, the errors measured against the reference solutions exhibit approximately second-order convergence in the $L^2$ norm and first-order convergence in the $H^1$ norm. The reported errors include geometric approximation and velocity-interpolation effects in addition to finite element discretization errors.

	\begin{table}[H]
		\centering\small\setlength{\tabcolsep}{3pt}
		\caption{2D spatial convergence of the nominal solution and interface sensitivities on $K_2$, with the $n=256$ finite element solution as the reference.}
		\label{tab:convergence}
		\begin{tabular}{cccccccc}
			\toprule
			$n$ & Elements & Free DOFs & $h_{\max}$
			& $e_U^{L^2}$ & $e_U^{H^1}$ & $e_S^{L^2}$ & $e_S^{H^1}$\\
			\midrule
			8 & 128 & 499 & 0.9442
			& $2.114\times10^{-2}$ & $8.791\times10^{-2}$
			& $2.499\times10^{-3}$ & $2.397\times10^{-2}$\\
			16 & 512 & 2147 & 0.4721
			& $5.654\times10^{-3}$ & $4.987\times10^{-2}$
			& $6.433\times10^{-4}$ & $1.056\times10^{-2}$\\
			32 & 2048 & 8899 & 0.2361
			& $1.319\times10^{-3}$ & $2.817\times10^{-2}$
			& $1.621\times10^{-4}$ & $5.326\times10^{-3}$\\
			64 & 8192 & 36227 & 0.1180
			& $3.254\times10^{-4}$ & $1.311\times10^{-2}$
			& $3.951\times10^{-5}$ & $2.623\times10^{-3}$\\
			128 & 32768 & 146179 & 0.0590
			& $6.669\times10^{-5}$ & $5.859\times10^{-3}$
			& $8.218\times10^{-6}$ & $1.175\times10^{-3}$\\
			\bottomrule
		\end{tabular}
	\end{table}

	\begin{table}[H]

		\centering\small\setlength{\tabcolsep}{3pt}
		\caption{3D spatial convergence of the nominal solution and interface sensitivities on $K_3$, with the $n=32$ finite element solution as the reference.}
		\label{tab:convergence3d}
		\begin{tabular}{cccccccc}
			\toprule
			$n$ & Elements & Free DOFs & $h_{\max}$ & $e_U^{L^2}$ & $e_U^{H^1}$ & $e_S^{L^2}$ & $e_S^{H^1}$\\
			\midrule
			4 &   384 &  1056 & 0.866 & 0.1154 & 0.4355 & 0.01986 & 0.08329\\
			8 &  3072 & 10468 & 0.433 & 0.02735 & 0.2384 & \ensuremath{6.186\times10^{-3}} & 0.04320 \\
			12 & 10368 & 37832 & 0.289 & 0.01126 & 0.1727 & \ensuremath{2.780\times10^{-3}} & 0.02997 \\
			16 & 24576 & 92748 & 0.217 & \ensuremath{5.765\times10^{-3}} & 0.1011 & \ensuremath{1.185\times10^{-3}} & 0.01835 \\
			\bottomrule
		\end{tabular}
	\end{table}


	\subsection{Perturbation order with exact discrete probability sums}\label{sec:exp_epsilon}
	
	This subsection investigates the accuracy of the shape expansion as $\epsilon\to0$ on a fixed mesh. We use exact probability summation for two discrete probability laws to assess the convergence rates and examine symmetry-induced cancellation.

	For the circular interface, we choose $\psi_1(\varphi)=0.35$ and $\psi_2(\varphi)=0.35\cos\varphi$. For the planar interface, we set
	\[
	\begin{aligned}
		\psi_1&=0.35a(x_1)a(x_2)
		\sin\!\left(\frac{\pi(x_1+1)}2\right)
		\sin\!\left(\frac{\pi(x_2+1)}2\right),\\
		\psi_2&=0.35a(x_1)a(x_2)
		\sin\!\left(\pi(x_1+1)\right)
		\sin\!\left(\frac{\pi(x_2+1)}2\right).
	\end{aligned}
	\]
	The random perturbation is given by $\beta=\xi_1\psi_1+\xi_2\psi_2$, where $\xi_1$ and $\xi_2$ are independent and identically distributed under either of two probability laws. The symmetric law assigns equal probabilities to $\xi=\pm1$, while the centered asymmetric law is defined by
	\[
	\mathbb P(\xi=-1/\sqrt2)=2/3,\quad
	\mathbb P(\xi=\sqrt2)=1/3.
	\]
	Both laws have zero mean and unit variance, with third moments $0$ and $1/\sqrt{2}$, respectively. Thus, they yield the same leading-order covariance, while higher-order terms may differ depending on the input law and response. Independence gives four joint realizations for each law, allowing exact computation of the discrete nonlinear moments without sampling error.

	The main experiments use meshes with $n=24$ in 2D and $n=6$ in 3D, together with the six perturbation amplitudes listed in Table~\ref{tab:epsilon_rates}. Tensor-product quadrature of orders 8 and 5 is used on each observation box in 2D and 3D, respectively. Mean errors are measured relative to $U_{0,h}$, and displacement variance is defined as the sum of the componentwise variances. Covariance errors are evaluated using a quadrature-based $L^2(K\times K)$ Frobenius norm, with the $up$ error restricted to the displacement--pressure covariance block. Since the coefficient vector has the identity covariance matrix under both laws, the leading-order covariance approximation is $\epsilon^2S_hS_h^*$. The reported rates describe the errors in the discrete perturbation expansion as measured by this quadrature on the observation set.

	\begin{table}[H]

		\centering\small\setlength{\tabcolsep}{3pt}
		\caption{Observed log--log slopes fitted over $\epsilon\in\{0.0025,0.005,0.01,0.02,0.04,0.06\}$. Each nonlinear reference statistic is computed by summing over all four states on the same finite element mesh. $C_{up}$ denotes the displacement--pressure covariance block.}
		\label{tab:epsilon_rates}
		\begin{tabular}{rccccccc}
			\toprule
			$d$ & Law & Mean $u$ & Mean $p$ & Var.\ $u$ & Var.\ $p$ & Full $C$ & $C_{up}$\\
			\midrule
			2 & asymmetric & 2.003 & 2.002 & 3.034 & 2.983 & 2.988 & 3.016\\
			2 & symmetric & 2.000 & 2.000 & 4.001 & 3.996 & 3.998 & 4.001\\
			3 & asymmetric & 2.000 & 1.998 & 2.998 & 2.987 & 2.989 & 2.996\\
			3 & symmetric & 2.000 & 2.000 & 3.999 & 4.000 & 4.000 & 4.000\\
			\bottomrule
		\end{tabular}
	\end{table}

\begin{figure}[tbp]

		\centering
		\includegraphics[width=0.94\linewidth]{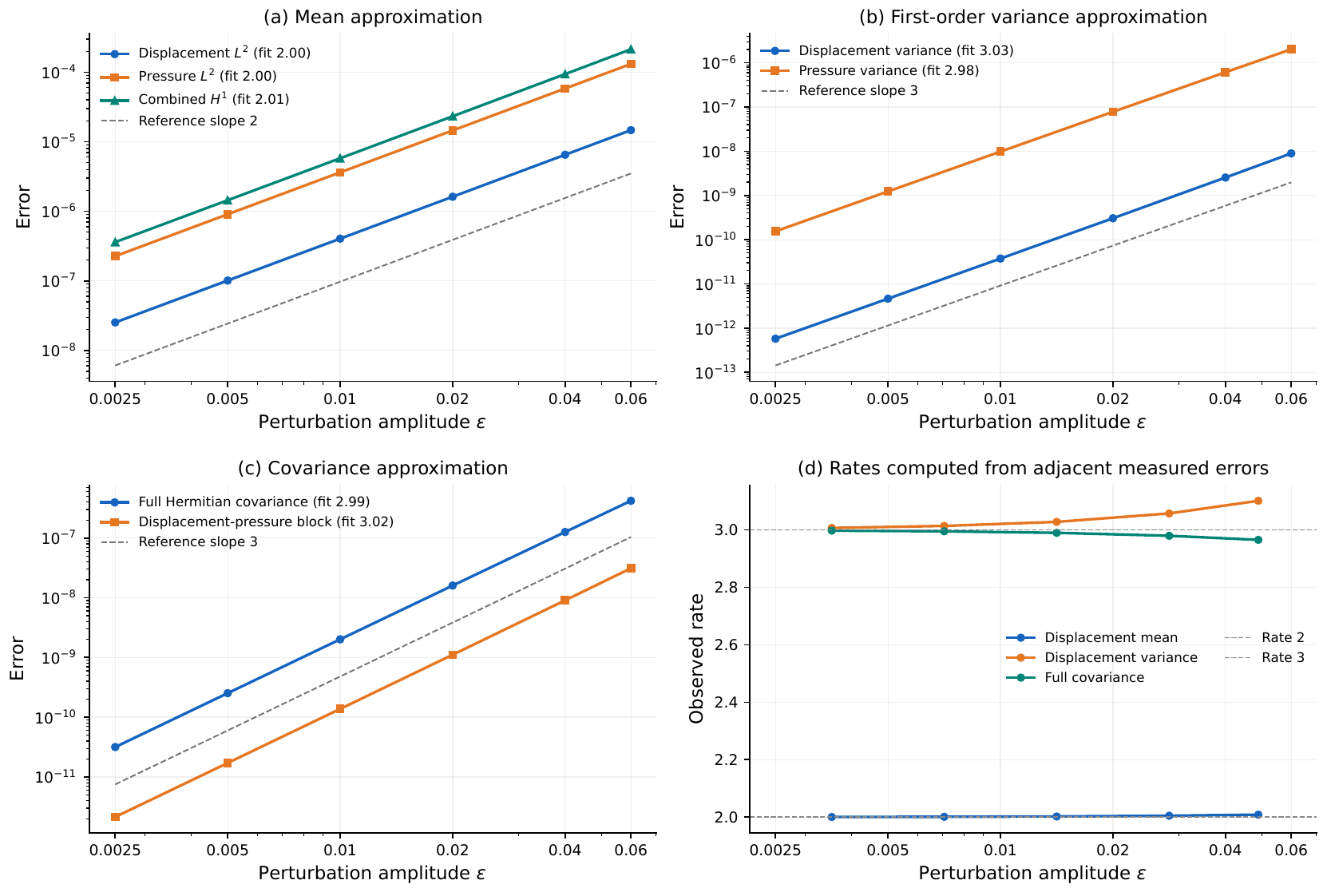}
		\caption{2D perturbation test under the asymmetric law on a fixed mesh ($n=24$): (a) mean approximation errors; (b) leading-order variance approximation errors; (c) covariance approximation errors; and (d) observed rates computed from adjacent perturbation amplitudes.}
		\label{fig:epsilon2d}
	\end{figure}

\begin{figure}[tbp]

		\centering
		\includegraphics[width=0.9\linewidth]{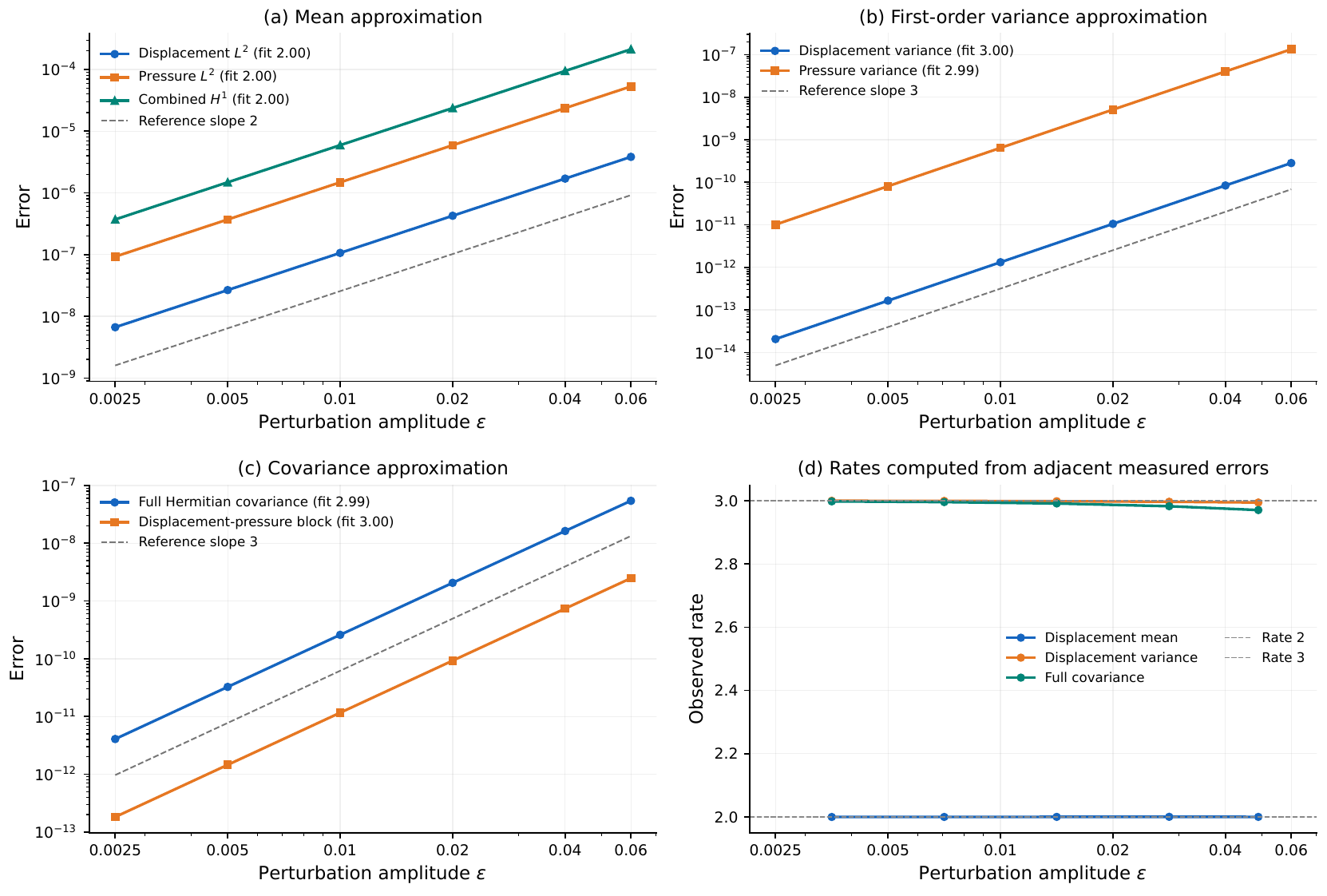}
		\caption{3D perturbation test under the asymmetric law on a fixed mesh ($n=6$): (a) mean approximation errors; (b) leading-order variance approximation errors; (c) covariance approximation errors; and (d) observed rates computed from adjacent perturbation amplitudes.}
		\label{fig:epsilon3d}
	\end{figure}

	For the asymmetric law, the mean errors exhibit approximately second-order convergence, while the variance and covariance errors converge at third order in both dimensions, consistent with the theoretical estimates. 
	Figures~\ref{fig:epsilon2d} and~\ref{fig:epsilon3d} show the log--log error plots and the slopes between successive amplitudes for the asymmetric law.
	For the symmetric law, the observed variance and covariance errors exhibit approximately fourth-order convergence, consistent with the cancellation of the cubic terms under the additional regularity assumptions in Remark~\ref{rem:symmetric_inputs}.


	Repeating the tests for $\epsilon=0.01,0.02,0.04$ on additional meshes with $n=16,48$ in 2D and $n=4$ in 3D yields similar observed convergence orders of 2, 3, and 4, with mesh-dependent constants. For example, at $\epsilon=0.04$, the 2D pressure-variance errors under the asymmetric law are $6.63\times10^{-7}$, $6.12\times10^{-7}$, and $5.84\times10^{-7}$ for $n=16,24,48$, respectively. These Taylor errors are measured at fixed mesh size. Table~\ref{tab:convergence} reports the separate spatial errors for the nominal solution and its sensitivities in 2D.

	\subsection{Low-rank approximation and coupled statistics}
	\label{sec:exp_lowrank}
	
	This subsection evaluates the errors introduced by low-rank approximation of the input covariance, using the full-rank leading-order statistics as a reference. It also examines the effect of neglecting the displacement--pressure cross-covariance in the joint output covariance.

	The variance error is measured in the $L^2(K)$ norm of the componentwise variance vector $(\operatorname{Var}u_1,\ldots,\operatorname{Var}u_d,\operatorname{Var}p)$, whereas the displacement variance in Table~\ref{tab:epsilon_rates} is the scalar displacement variance $V_{u,\epsilon}$ in \eqref{eq:variance_direct}, which sums the displacement component variances. To assess low-rank approximation errors, we consider $q=17$ circular-interface modes,
	\[
	\boldsymbol\psi(\phi)
	=0.35(1,\cos\phi,\sin\phi,\ldots,\cos8\phi,\sin8\phi)^T.
	\]
	Let $\nu=(0,1,1,\ldots,8,8)$ and set
	\[
	d_a=(1+\nu_a^2)^{-3/2},\quad
	R_{ab}=0.65^{|a-b|},\quad
	Q=C_\xi=
	\frac{\operatorname{diag}(d)R\operatorname{diag}(d)}
	{\operatorname{tr}(\operatorname{diag}(d)R\operatorname{diag}(d))}.
	\]
	The covariance is realized by $\boldsymbol\xi=Q^{1/2}\boldsymbol\zeta$, where the components of $\boldsymbol\zeta$ are independent and uniformly distributed on $[-\sqrt{3},\sqrt{3}]$. For the leading-order compression, we apply pivoted Cholesky factorization directly to $Q$, with residual trace $\eta_m=\operatorname{tr}(Q-L_mL_m^T)$ as defined in \eqref{eq:coefficient_cholesky_residual}. For each rank $m$, we solve $m$ sensitivity problems and compare the resulting statistics with the full-rank reference using all $q$ sensitivity directions.

	\begin{table}[H]

		\centering\small\setlength{\tabcolsep}{3pt}
		\caption{Seventeen-mode covariance compression for the 2D model with $n=24$. Relative errors are measured against full-rank leading-order statistics, with $\epsilon^2$ removed. Timings include nominal and sensitivity solves and statistical evaluation.}
		\label{tab:lowrank}
		\begin{tabular}{ccccc}
			\toprule
			$m$ & $\eta_m$ & Rel.\ variance error & Rel.\ covariance error & Measured time (s)\\
			\midrule
			1 & 0.152 & \ensuremath{4.54\times10^{-3}} & \ensuremath{4.09\times10^{-3}} & 0.063\\
			2 & 0.0508 & \ensuremath{1.44\times10^{-3}} & \ensuremath{1.89\times10^{-3}} & 0.072\\
			4 & \ensuremath{4.10\times10^{-3}} & \ensuremath{1.47\times10^{-5}} & \ensuremath{2.61\times10^{-5}} & 0.103\\
			8 & \ensuremath{1.76\times10^{-4}} & \ensuremath{3.95\times10^{-8}} & \ensuremath{4.47\times10^{-8}} & 0.129\\
			12 & \ensuremath{2.06\times10^{-5}} & \ensuremath{4.85\times10^{-10}} & \ensuremath{4.77\times10^{-10}} & 0.165\\
			17 & 0 & \ensuremath{1.45\times10^{-14}} & \ensuremath{1.85\times10^{-14}} & 0.220\\
			\bottomrule
		\end{tabular}
	\end{table}

	\begin{figure}[H]

		\centering
		\includegraphics[width=0.94\linewidth]{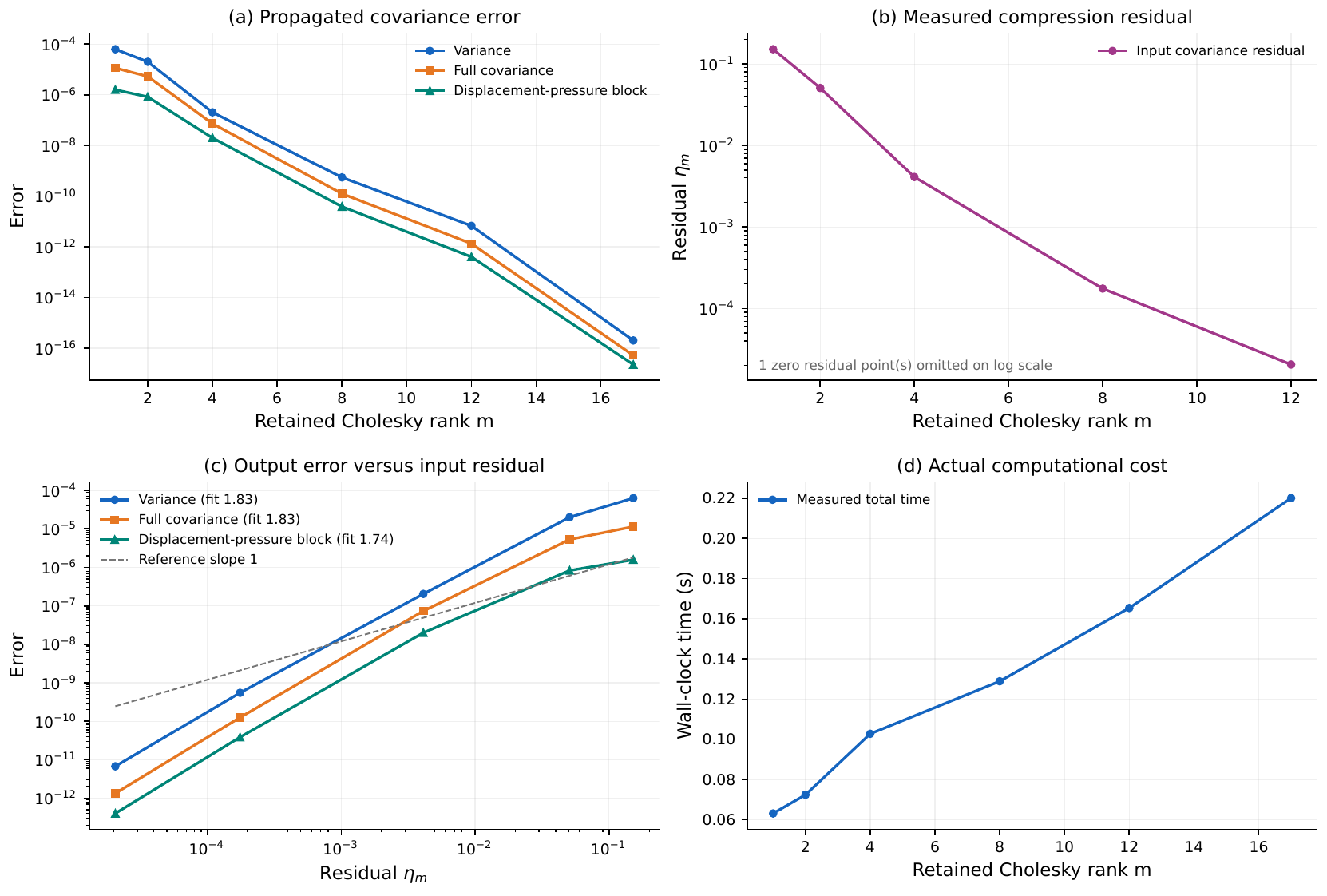}
		\caption{Rank dependence of covariance compression for the 17-mode input. Output errors are measured against full-rank leading-order statistics, with the input covariance residual shown separately.}
		\label{fig:lowrank}
	\end{figure}

	Table~\ref{tab:lowrank} and Figure~\ref{fig:lowrank} show that the output errors decrease rapidly as the retained rank increases. At $m=8$, the residual trace is $1.76\times10^{-4}$, while the relative variance and covariance errors are $3.95\times10^{-8}$ and $4.47\times10^{-8}$, respectively. At $m=12$, both errors fall below $5\times10^{-10}$. The displacement--pressure covariance block exhibits similar error decay. For the five ranks with positive residual trace ($m=1,2,4,8,12$), log--log fits of the relative errors against $\eta_m$ give slopes of approximately $1.83$ for variance and full covariance and $1.74$ for the cross-covariance block. Meanwhile, the measured time increases from $0.129$ s at $m=8$ to $0.220$ s at full rank.
	At $m=17$, the relative errors reach approximately $10^{-14}$. The compressed sensitivity solutions also agree with the corresponding linear combinations of the full-rank sensitivities to within $10^{-14}$ in the measured norm. These results demonstrate the accuracy and computational efficiency of covariance compression for the prescribed covariance structure and observation set.

	To illustrate the effect of displacement--pressure cross-covariance, consider the spatial averages $Q_u=|K|^{-1}\int_K u_1(\mathbf{x}){\mathrm{d}\,\mathbf{x}}$ and $Q_p=|K|^{-1}\int_K p(\mathbf{x}){\mathrm{d}\,\mathbf{x}}$. For their leading-order centered sensitivities $q_u$ and $q_p$, define
	\[
	\rho_{up}=
	\frac{\mathbb E[q_u\overline{q_p}]}
	{\sqrt{\mathbb E|q_u|^2\,\mathbb E|q_p|^2}},
	\quad
	Q_{\rm joint}=
	\frac{q_u}{\sqrt{\mathbb E|q_u|^2}}
	+\frac{q_p}{\sqrt{\mathbb E|q_p|^2}}.
	\]
	Then $\operatorname{Var}(Q_{\rm joint}) =2+2\operatorname{Re}\rho_{up}$. We use the two perturbation directions from Section~\ref{sec:exp_epsilon} without phase adjustment. Omitting the cross-covariance gives a joint variance of $2$, which exceeds the coupled value by $19.4\%$ in 2D and $95.2\%$ in 3D (Table~\ref{tab:coupled_qoi}). These results demonstrate the importance of retaining displacement--pressure cross-covariance in joint statistical estimates.

	\begin{table}[H]

		\centering\small\setlength{\tabcolsep}{3pt}
		\caption{Displacement--pressure correlation and normalized joint-output variance. Neglecting cross-covariance gives a variance of $2$.}
		\label{tab:coupled_qoi}
		\begin{tabular}{cccccc}
			\toprule
			$d$ & $\operatorname{Re}\rho_{up}$ & $\operatorname{Im}\rho_{up}$ & $|\rho_{up}|$ & Joint variance & Overestimate\\
			\midrule
			2 & -0.1622 & -0.4128 & 0.4435 & 1.6756 & 19.4\%\\
			3 & -0.4876 & -0.6093 & 0.7804 & 1.0248 & 95.2\%\\
			\bottomrule
		\end{tabular}
	\end{table}

	\subsection{Computational efficiency and sampling uncertainty}
	\label{sec:exp_cost}
	
	This subsection compares statistical accuracy and computational cost on a fixed spatial mesh and evaluates the sampling uncertainty of Monte Carlo estimates. The study focuses on probability approximation and shape-expansion errors, excluding the common spatial discretization error.

	We use the 2D mesh with $n=16$, $\epsilon=0.04$, and the first three circular-interface modes. The coefficient covariance is constructed from the formula in Section~\ref{sec:exp_lowrank} using $q=3$ and frequencies $(0,1,1)$, with its trace renormalized to one; it is not the unrenormalized principal block of the $q=17$ matrix. Writing $C_\xi=LL^T$, we generate $\boldsymbol\xi=L\boldsymbol\zeta$, where the components of $\boldsymbol\zeta$ are independent and uniformly distributed on $[-\sqrt3,\sqrt3]$. We compare analytic sensitivities, centered finite-difference sensitivities with step size $10^{-4}$, tensor-product Gauss quadrature with $2,3,4,6,$ or $8$ points per coordinate, and Monte Carlo sampling with $N=16,\ldots,1024$. The Monte Carlo samples are generated using NumPy's default random generator with seed $20261021$, and sample covariances use the unbiased estimator with denominator $N-1$.

	The reference statistics are computed on the same mesh using tensor-product Gauss quadrature with eight points per coordinate (512 nonlinear solves). The relative differences between the $6^3$- and $8^3$-point rules are $2.00\times10^{-14}$, \ensuremath{1.36\times10^{-11}}, and \ensuremath{1.26\times10^{-11}} for the mean, variance, and full covariance, respectively, confirming the accuracy of the reference probability quadrature.

	The reported timings include PDE matrix and right-hand-side assembly, factorization, and solves, with the nominal factorization reused for sensitivity computations. Common mesh and finite element space setup and statistical postprocessing are excluded.

	\begin{table}[H]

		\centering\small\setlength{\tabcolsep}{3pt}
		\caption{Statistical accuracy and computational cost for the 2D model ($n=16$, $q=3$, $\epsilon=0.04$). Relative errors are measured against the $8^3$-point tensor Gauss reference statistics.}
		\label{tab:cost}
		\begin{tabular}{lrrrrr}
			\toprule
			Method & Factorizations & Mean error & Variance error & Covariance error & Time (s)\\
			\midrule
			Analytic sensitivity & 1 & \ensuremath{2.55\times10^{-4}} & \ensuremath{3.43\times10^{-4}} & \ensuremath{8.61\times10^{-4}} & 0.0268\\
			Central-difference sensitivity & 7 & \ensuremath{2.55\times10^{-4}} & \ensuremath{3.43\times10^{-4}} & \ensuremath{8.61\times10^{-4}} & 0.1359\\
			Tensor Gauss, $2^3$ & 8 & \ensuremath{9.83\times10^{-8}} & \ensuremath{3.74\times10^{-5}} & \ensuremath{3.89\times10^{-4}} & 0.1346\\
			Tensor Gauss, $3^3$ & 27 & \ensuremath{1.04\times10^{-9}} & \ensuremath{7.14\times10^{-7}} & \ensuremath{6.96\times10^{-7}} & 0.5618\\
			Tensor Gauss, $6^3$ & 216 & \ensuremath{2.00\times10^{-14}} & \ensuremath{1.36\times10^{-11}} & \ensuremath{1.26\times10^{-11}} & 3.8436\\
			Monte Carlo & 64 & \ensuremath{3.33\times10^{-3}} & 0.0342 & 0.0295 & 1.2605\\
			Monte Carlo & 256 & \ensuremath{6.95\times10^{-4}} & 0.0445 & 0.0439 & 4.6659\\
			Monte Carlo & 1024 & \ensuremath{4.73\times10^{-4}} & 0.0144 & 0.0145 & 17.6694\\
			\bottomrule
		\end{tabular}
	\end{table}
	
	\begin{figure}[H]

		\centering
		\includegraphics[width=0.94\linewidth]{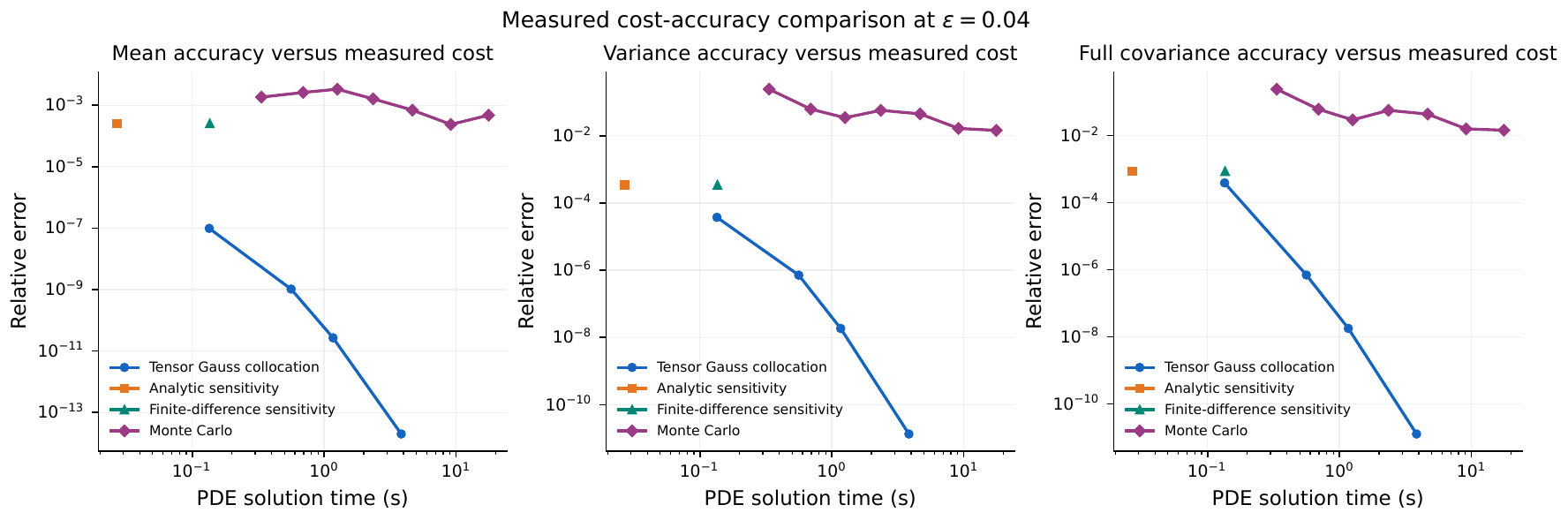}
		\caption{Relative errors in the mean, variance, and full covariance versus PDE solution time for the 2D model ($n=16$, $q=3$, $\epsilon=0.04$).}
		\label{fig:cost}
	\end{figure}

	Table~\ref{tab:cost} and Figure~\ref{fig:cost} compare the accuracy and computational cost of these methods. Analytic and finite-difference sensitivities yield nearly identical statistical errors, with a relative difference of $2.48\times10^{-9}$ between their sensitivity fields. The analytic computation used one factorization and took $0.0268$ s. Centered differences used seven factorizations and took $0.1359$ s, giving a measured time ratio of $0.1359/0.0268\approx5.07$. 
	Direct Gauss quadrature achieves higher accuracy at increased computational cost. For example, the $3^3$-point rule reduces the relative mean, variance, and covariance errors to $1.04\times10^{-9}$, $7.14\times10^{-7}$, and $6.96\times10^{-7}$, respectively, with a computational time of $0.5618$ s. In contrast, Monte Carlo errors fluctuate with sample size and remain substantially larger. At $N=1024$, Monte Carlo takes 17.669 s, with relative mean, variance, and covariance errors of $4.73\times10^{-4}$, $0.0144$, and $0.0145$, respectively. These results illustrate the efficiency of analytic sensitivities for the prescribed accuracy, while direct quadrature provides higher accuracy in this low-dimensional setting.

	\begin{table}[H]

		\centering\small\setlength{\tabcolsep}{3pt}
		\caption{Sampling variability of the seeded Monte Carlo stream. RMS values are plug-in sampling-error estimates rather than confidence intervals.}
		\label{tab:mc_uncertainty}
		\begin{tabular}{ccccc}
			\toprule
			$N$ & Observed mean error & Mean RMS estimate & Observed variance error & Variance RMS estimate\\
			\midrule
			64 & \ensuremath{3.33\times10^{-3}} & \ensuremath{1.56\times10^{-3}} & 0.0342 & 0.13\\
			256 & \ensuremath{6.95\times10^{-4}} & \ensuremath{7.89\times10^{-4}} & 0.0445 & 0.0586\\
			1024 & \ensuremath{4.73\times10^{-4}} & \ensuremath{3.83\times10^{-4}} & 0.0144 & 0.0283\\
			\bottomrule
		\end{tabular}
	\end{table}

	Table~\ref{tab:mc_uncertainty} reports plug-in RMS sampling-error estimates based on the sample variance for the mean and the fourth central moment for the variance. The estimated errors decrease approximately as $N^{-1/2}$, while the observed variance error increases from $N=64$ to $N=256$ and then decreases at $N=1024$. At $N=1024$, the relative RMS estimates for the mean and variance are $3.83\times10^{-4}$ and $2.83\times10^{-2}$, respectively. These estimates characterize sampling variability rather than total approximation errors. Additional variability summaries from disjoint batches with $N=16,32,64,128$, together with per-solve logs and independent timing results, are provided with the accompanying data.

	\section{Conclusion}\label{section:Conclusion}

	{ We developed a moment analysis for a reduced time-harmonic Biot system with random material interfaces. The} model retains mixture inertia, quasi-static Darcy discharge, and explicit storage and coupling coefficients{ . The shape derivative was stated in a total-traction form that explicitly accounts for normal variation, coefficient jumps, endpoint compatibility, and} source terms. {  Under a first-order shape Taylor expansion, the nominal solution approximates the mean with $O(\epsilon^2)$ accuracy,} while $\epsilon^2$ times the Hermitian covariance of the first-order shape derivative approximates the solution covariance with $O(\epsilon^3)$ accuracy. The corresponding direct variance approximations achieve the same accuracy. { Pivoted Cholesky compression of the interface covariance reduces the variance calculation to deterministic shape-derivative solves.} This work extends an established moment framework to the coupled Biot system, with a rigorous analysis of its shape sensitivities and statistical moments.

	\vspace{10mm}
	
	\noindent {\bf Acknowledgements.}
	The work of F. Huo is supported in part by the National Key R\&D Program of China (2024YFA1012302). The work of J. Li  was partially supported by the Shenzhen Sci-Tech Fund No. JCYJ20250604144225032 and  RCJC20200714114556020, Guangdong Basic and Applied Research Fund No. 2023B1515250005,  National Center for Applied Mathematics Shenzhen, and SUSTech International Center for Mathematics. The work of K. Zhang is supported in part by the National Natural Science Foundation of China (Grant No. 12271207).

	\vspace{5mm}
	
	\noindent {\bf Declaration of competing interest.}  The authors have no relevant financial or non-financial interests to disclose.

	\appendix
	
	\section{Fixed-domain operator formulation for the shape Taylor expansion}\label{app:fixed_domain}

	\begin{lemma}\label{lem:fredholm}
		Under the coefficient assumptions of Section~\ref{section:Model}, the operator $A_0:X\to X^*$ defined by $\langle A_0U,W\rangle=b(U,W)$ is Fredholm of index zero. Under Assumption~\ref{ass:data}, it is an isomorphism.
	\end{lemma}
	\begin{proof}
		The principal form
		\[
		a(U,W)=\sum_i\int_{D_i}
		\bigl(2\mu_i\boldsymbol\varepsilon(\mathbf u):
		\overline{\boldsymbol\varepsilon(\mathbf v)}
		+\lambda_i\nabla\cdot\mathbf u\,\overline{\nabla\cdot\mathbf v}
		+\kappa_i\nabla p\cdot\overline{\nabla q}\bigr)\,{d\mathbf x}
		\]
		is bounded and coercive on $X$. Indeed, the conditions $\mu_i>0$ and $d\lambda_i+2\mu_i>0$ give a uniform positive lower bound for the elastic energy on symmetric strains, by decomposition into trace-free and spherical parts. This bound, together with $\min_i\kappa_i>0$ and Korn's and Poincar\'e's inequalities, gives a constant $c_a>0$ such that
		\[
		a(U,U)\ge c_a\bigl(\|\mathbf u\|_{H^1(D)}^2+\|p\|_{H^1(D)}^2\bigr)
		=c_a\|U\|_X^2 .
		\]
		The Lax--Milgram theorem then gives an isomorphism from $X$ onto $X^*$. The terms with $\rho_i$ and $c_i$ are compact by the compact embedding $H^1_0(D)\hookrightarrow L^2(D)$. For the first coupling, the map $p\mapsto\alpha p$ is compact from $H^1_0(D)$ into $L^2(D)$. For $h\in L^2(D)$, the functional $\mathbf v\mapsto-\int_D h\,\overline{\nabla\cdot\mathbf v}\,d\mathbf x$ satisfies
		\[
		\left|\int_D h\,\overline{\nabla\cdot\mathbf v}\,d\mathbf x\right|
		\le C\|h\|_{L^2(D)}\|\mathbf v\|_{H^1(D)}.
		\]
		Composing the resulting bounded map into the displacement anti-dual with $p\mapsto\alpha p$ gives the compact coupling operator. 
		The second coupling is the composition of the bounded map $\mathbf u\mapsto\alpha\nabla\cdot\mathbf u$ into $L^2(D)$ and the compact inclusion $L^2(D)\hookrightarrow H^{-1}(D)$, up to the constant factor $i\theta$. Thus $A_0$ is a compact perturbation of an isomorphism. Fredholm index zero and the assumed injectivity prove the assertion.
	\end{proof}

	\begin{lemma}
		\label{lem:low_frequency_nonresonance}
		Suppose $\theta>0$ and $\lambda_i\ge0$. Set $\mu_{\min}=\min_i\mu_i$ and $\rho_{\max}=\max_i\rho_i$. Let $\lambda_1(D)$ be the first scalar Dirichlet eigenvalue of $-\Delta$ on $D$. If
		\begin{equation}\label{eq:low_frequency_condition}
			\theta^2\rho_{\max}<\mu_{\min}\lambda_1(D),
		\end{equation}
		the injectivity condition in Assumption~\ref{ass:data} holds. For a rectangular box,
		\[
		\lambda_1(D)=\pi^2\sum_{j=1}^dL_j^{-2},
		\]
		where $L_j$ are its side lengths.
	\end{lemma}
	\begin{proof}
		For a homogeneous solution, use $W=(\mathbf u,p/(i\theta))$ in \eqref{eq:weak_state}. Since the form is conjugate-linear in its test argument and the coefficients $\alpha_i$ are real, the two coupling terms sum to
		\[
		-\sum_i\int_{D_i}\alpha_i
		\bigl(p\,\overline{\nabla\cdot\mathbf u}
		+(\nabla\cdot\mathbf u)\,\overline p\bigr)\,d\mathbf x,
		\]
		which is real. The elastic, inertial, and storage terms are also real.
		The imaginary part is
		\[
		\theta^{-1}\sum_i\int_{D_i}\kappa_i|\nabla p|^2\,{d\mathbf x}=0.
		\]
		Since each $\kappa_i$ is positive, $\nabla p=0$. The zero boundary trace of $p$ then gives $p=0$. The remaining identity is
		\[
		a_{\rm el}(\mathbf u,\mathbf u)
		-\theta^2\sum_i\int_{D_i}\rho_i|\mathbf u|^2\,{d\mathbf x}=0,
		\]
		where $a_{\rm el}$ is the elastic part of $a$. The global identity $2\|\boldsymbol\varepsilon(\mathbf u)\|_{L^2}^2 =\|\nabla\mathbf u\|_{L^2}^2+\|\nabla\cdot\mathbf u\|_{L^2}^2$ for $\mathbf u\in H^1_0(D;\mathbb C^d)$, together with the scalar Dirichlet Poincar\'e inequality applied to each displacement component, gives
		\[
		a_{\rm el}(\mathbf u,\mathbf u)
		\ge 2\mu_{\min}\|\boldsymbol\varepsilon(\mathbf u)\|_{L^2}^2
		\ge\mu_{\min}\|\nabla\mathbf u\|_{L^2}^2
		\ge\mu_{\min}\lambda_1(D)\|\mathbf u\|_{L^2}^2.
		\]
		Here the first inequality uses $\lambda_i\ge0$. Combining this estimate with the displacement identity gives
		\[
		\bigl(\mu_{\min}\lambda_1(D)-\theta^2\rho_{\max}\bigr)
		\|\mathbf u\|_{L^2(D)}^2\le0 .
		\]
		Condition~\eqref{eq:low_frequency_condition} implies $\mathbf u=0$.
	\end{proof}

	For a fixed $\beta$, write
	\[
	F_\epsilon=F_{\epsilon,\beta},\quad
	G_\epsilon=(DF_\epsilon)^{-1},\quad J_\epsilon=\det DF_\epsilon .
	\]
	For a vector field $\mathbf v$ and a scalar field $q$, the transformed differential operators are
	\[
	\begin{aligned}
		\nabla_\epsilon\mathbf v&=\nabla\mathbf v\,G_\epsilon,
		&\operatorname{div}_\epsilon\mathbf v
		&=\operatorname{tr}(\nabla\mathbf v\,G_\epsilon),\\
		\boldsymbol\varepsilon_\epsilon(\mathbf v)
		&=\operatorname{sym}(\nabla\mathbf v\,G_\epsilon),
		&\nabla_\epsilon q&=G_\epsilon^T\nabla q .
	\end{aligned}
	\]
	The pullback form $b_\epsilon$ is obtained from \eqref{eq:bilinear_form} by replacing the gradients, divergences, and symmetric gradients with these transformed quantities and multiplying the integrand by $J_\epsilon$. The material coefficients remain fixed on each reference subdomain $D_i$.
	The transformed load is
	\begin{equation}\label{eq:pullback_load}
		\ell_\epsilon(W)=\int_D J_\epsilon
		\bigl((\mathbf f\circ F_\epsilon)\cdot\overline{\mathbf v}
		+(g\circ F_\epsilon)\overline q\bigr)\,d\mathbf x.
	\end{equation}
	Thus
	\begin{equation}\label{eq:pullback_weak}
		A_\epsilon\widetilde U_\epsilon=\ell_\epsilon\quad\text{in }X^*,
		\quad \langle A_\epsilon U,W\rangle=b_\epsilon(U,W).
	\end{equation}
	All terms in \eqref{eq:pullback_weak} act on the fixed space $X$ through volume integrals.

	\begin{theorem}\label{thm:pullback}
		Under Assumption~\ref{ass:data} and the geometric, extension, and measurability assumptions of Section~\ref{sec:stochastic_geometry}, there exists $\epsilon_0>0$ such that, for every $\beta\in B$ with $\|\beta\|_B\le M_\beta$ and every $|\epsilon|\le\epsilon_0$, the operator $A_\epsilon:X\to X^*$ is invertible and satisfies
		\[
		\|A_\epsilon^{-1}\|_{\mathcal L(X^*,X)}
		\le2\|A_0^{-1}\|_{\mathcal L(X^*,X)}.
		\]
		By~\eqref{eq:beta_assumption}, the same $\epsilon_0$ applies to $\beta(\omega)$ for almost every $\omega\in\Xi$. The map $\epsilon\mapsto\widetilde U_\epsilon$ is $C^2$ in $X$. Its material derivative solves
		\begin{equation}\label{eq:material_derivative}
			A_0\dot U[\beta]=\ell'_0[\beta]-A'_0[\beta]U_0.
		\end{equation}
		The map $\beta\mapsto\dot U[\beta]$ is real-linear and continuous from $B$ to $X$, and
		\begin{equation}\label{eq:material_taylor}
			\|\dot U[\beta]\|_X\le C\|\beta\|_B,\quad
			\|\widetilde U_\epsilon-U_0-\epsilon\dot U[\beta]\|_X
			\le C\epsilon^2\|\beta\|_B^2.
		\end{equation}
		The constants depend on $M_\beta$, the nominal resolvent, the coefficients, the data, and the extension operator. They are uniform over the stated ranges of $\beta$ and $\epsilon$.
	\end{theorem}
	\begin{proof}
		The extension estimate $\|V_\beta\|_{W^{2,\infty}(D)}\le C_E\|\beta\|_B$ gives uniform control of $V_\beta$ and $DV_\beta$ on the ball $\|\beta\|_B\le M_\beta$. For sufficiently small $|\epsilon|$, the Jacobian matrix $DF_\epsilon=I+\epsilon DV_\beta$ and its inverse are uniformly bounded. Smooth dependence of the inverse and determinant on $DF_\epsilon$ therefore gives $C^2$ dependence of the pulled-back form coefficients on $\epsilon$ in $L^\infty$. The fixed $C^2$ sources likewise give a $C^2$ load in $X^*$. In particular,
		\[
		\|A_\epsilon-A_0\|_{\mathcal L(X,X^*)}
		\le C|\epsilon|\|\beta\|_B.
		\]
		Choose $\epsilon_{\rm geo}>0$ such that, for $|\epsilon|<\epsilon_{\rm geo}$, the diffeomorphism property and the bounds on $DF_\epsilon$ and $G_\epsilon$ hold uniformly. By Lemma~\ref{lem:fredholm}, $A_0$ is invertible. 
		Choose $0<\epsilon_0<\epsilon_{\rm geo}$ so that $C\epsilon_0 M_\beta\|A_0^{-1}\|_{\mathcal L(X^*,X)}\le 1/2$. Then $\|A_0^{-1}(A_\epsilon-A_0)\|_{\mathcal L(X,X)}\le 1/2$ throughout the stated ranges of $\epsilon$ and $\beta$. The Neumann series gives
		\[
		A_\epsilon^{-1}
		=\bigl[I+A_0^{-1}(A_\epsilon-A_0)\bigr]^{-1}A_0^{-1},
		\quad
		\|A_\epsilon^{-1}\|_{\mathcal L(X^*,X)}
		\le 2\|A_0^{-1}\|_{\mathcal L(X^*,X)}.
		\]
		Since inversion is smooth on the set of isomorphisms, the map $\epsilon\mapsto A_\epsilon^{-1}\ell_\epsilon=\widetilde U_\epsilon$ is $C^2$ in $X$. Differentiating~\eqref{eq:pullback_weak} at zero gives \eqref{eq:material_derivative}.
		At $\epsilon=0$, the first derivatives of the operator and load are linear in $V_\beta$, and their second derivatives are quadratic in $V_\beta$. Along the perturbation path, write $H=DV_\beta$. Then $G_\epsilon'=-G_\epsilon H G_\epsilon$ and $G_\epsilon''=2G_\epsilon H G_\epsilon H G_\epsilon$. Jacobi's formula gives
		\[
		\begin{aligned}
			J_\epsilon'&=J_\epsilon\operatorname{tr}(G_\epsilon H),\\
			J_\epsilon''&=J_\epsilon\left[
			\bigl(\operatorname{tr}(G_\epsilon H)\bigr)^2
			-\operatorname{tr}(G_\epsilon H G_\epsilon H)\right].
		\end{aligned}
		\]
		These identities, the uniform geometric bounds, and the chain rule for the fixed $C^2$ sources in Assumption~\ref{ass:data} yield
		\[
		\|\partial_\epsilon^j A_\epsilon\|_{\mathcal L(X,X^*)}
		+\|\partial_\epsilon^j\ell_\epsilon\|_{X^*}
		\le C\|\beta\|_B^j,\quad j=1,2,\quad |\epsilon|\le\epsilon_0.
		\]
		Here $C$ is uniform for $\|\beta\|_B\le M_\beta$.
		Differentiating~\eqref{eq:pullback_weak} along the path gives
		\[
		A_\epsilon\partial_\epsilon\widetilde U_\epsilon
		=\ell'_\epsilon-A'_\epsilon\widetilde U_\epsilon.
		\]
		A second differentiation gives
		\[
		A_\epsilon\partial_{\epsilon\epsilon}\widetilde U_\epsilon
		=\ell''_\epsilon-A''_\epsilon\widetilde U_\epsilon
		-2A'_\epsilon\partial_\epsilon\widetilde U_\epsilon.
		\]
		The pulled-back loads are uniformly bounded in $X^*$. Since $\widetilde U_\epsilon=A_\epsilon^{-1}\ell_\epsilon$, the states are uniformly bounded in $X$. Applying the inverse bound first to the equation for the first derivative and then to the equation for the second derivative yields $\|\partial_\epsilon\widetilde U_\epsilon\|_X\le C\|\beta\|_B$ and $\|\partial_{\epsilon\epsilon}\widetilde U_\epsilon\|_X \le C\|\beta\|_B^2$.
		Taylor's formula in $X$ gives the Bochner integral remainder
		\[
		\widetilde U_\epsilon-U_0-\epsilon\dot U[\beta]
		=\int_0^\epsilon(\epsilon-t)\,\partial_{tt}\widetilde U_t\,dt.
		\]
		The bound on the second derivative gives the remainder estimate in \eqref{eq:material_taylor} for both signs of $\epsilon$.
		Finally, $V_\beta=E\beta$ depends real-linearly on $\beta$, so \eqref{eq:material_derivative} gives real-linearity of $\beta\mapsto\dot U[\beta]$, and the first bound in \eqref{eq:material_taylor} gives continuity.
	\end{proof}

	The right-hand side of \eqref{eq:material_derivative} is explicit.
	In the formulas below, the argument $V=V_\beta$ denotes the extension of the interface amplitude $\beta$ in \eqref{eq:material_derivative}. With $H=DV$ and $\tau=\operatorname{tr}H$, define
	\[
	e'_V(\mathbf v)=-\operatorname{sym}(\nabla\mathbf v\,H),\quad
	d'_V(\mathbf v)=-\operatorname{tr}(\nabla\mathbf v\,H),\quad
	r'_V(q)=-H^T\nabla q.
	\]
	If $\mathcal B_i(U,W)$ denotes the integrand in braces in \eqref{eq:bilinear_form}, direct differentiation yields
	\begin{align} \label{eq:form_derivative}
		\begin{aligned}
			b'_0[V](U,W)=\sum_i\int_{D_i}\bigl\{&
			\tau\mathcal B_i(U,W)
			+2\mu_i\bigl(e'_V(\mathbf u):\overline{\boldsymbol\varepsilon(\mathbf v)}
			+\boldsymbol\varepsilon(\mathbf u):\overline{e'_V(\mathbf v)}\bigr)
			\\
			&+\lambda_i\bigl(d'_V(\mathbf u)\overline{\nabla\cdot\mathbf v}
			+(\nabla\cdot\mathbf u)\overline{d'_V(\mathbf v)}\bigr)
			-\alpha_i p\overline{d'_V(\mathbf v)}\\
			&+\kappa_i\bigl(r'_V(p)\cdot\overline{\nabla q}
			+\nabla p\cdot\overline{r'_V(q)}\bigr)
			+i\theta\alpha_i d'_V(\mathbf u)\overline q\bigr\}\,d\mathbf x.
		\end{aligned}
	\end{align}
	Differentiating the load in \eqref{eq:pullback_load} gives
	\begin{align}\label{eq:load_derivative}
		\ell'_0[V](W)=\int_D\bigl\{
		(\tau\mathbf f+D\mathbf f\,V)\cdot\overline{\mathbf v}
		+(\tau g+\nabla g\cdot V)\overline q\bigr\}\,d\mathbf x.
	\end{align}
	Even constant physical sources therefore give a generally nonzero Jacobian contribution to the load derivative.

\end{document}